\documentclass{amsart}
\usepackage{amsmath, amscd, amssymb}

\newtheorem{theorem}{Theorem}[section]
\newtheorem{lemma}[theorem]{Lemma}
\newtheorem{corollary}[theorem]{Corollary}

\newtheorem{proposition}[theorem]{Proposition}

\theoremstyle{definition}
\newtheorem{example}[theorem]{Example}
\newtheorem{remark}[theorem]{Remark}
\newtheorem{definition}[theorem]{Definition}

\def\L{\mathcal{L}}
\def\A{\mathcal{A}}
\def\U{\mathcal{U}}
\def\R{\mathcal{R}}
\def\S{\mathcal{S}}

\newcommand{\ZZ}{{\mathbb Z}}
\newcommand{\NN}{{\mathbb N}}
\newcommand{\FF}{{\mathbb F}}
\newcommand{\CC}{{\mathbb C}}

\def\max{\operatorname{max}}
\def\log{\operatorname{log}}
\def\ker{\operatorname{ker}}
\def\mod{\operatorname{mod}}
\def\rm{\operatorname{RM}}
\def\dm{\operatorname{dM}}
\def\th{\operatorname{Th}}
\def\alg{\operatorname{alg}}

\def\stab{\operatorname{Stab}}
\def\aut{\operatorname{Aut}}
\def\dim{\operatorname{dim}}
\def\deg{\operatorname{deg}}
\def\tp{\operatorname{tp}}
\def\ind{\operatorname{ind}}
\def\id{\operatorname{id}}

\def\c{\overline{c}}
\def\b{\overline{b}}
\def\a{\overline{a}}
\def\r{\overline{r}}
\def\s{\overline{s}}

\def\m{\overline{m}}
\def\y{\overline{y}}
\def\x{\overline{x}}
\def\z{\overline{z}}
\def\w{\overline{w}}

\begin{document}

\title[$F$-structures and Semiabelian Varieties]{$F$-structures and Integral Points on Semiabelian Varieties over Finite Fields}

\author{Rahim Moosa}
\address{The Fields Institute \\
222 College Street \\ 
Toronto, Ontario M5T 3J1 \\ 
Canada}
\address{University of California, Berkeley \\
Department of Mathematics \\
Evans Hall \\
Berkeley, California 94720-3840 \\
USA}
\curraddr{Massachusetts Institute of Technology \\
Department of Mathematics \\
77 Mass. Ave. \\
Cambridge, MA 02139-4307 \\
USA}
\email{moosa@math.mit.edu}

\author{Thomas Scanlon}
\address{University of California, Berkeley \\
Department of Mathematics \\
Evans Hall \\
Berkeley, California 94720-3840 \\
USA}
\email{scanlon@math.berkeley.edu}

\thanks{Rahim Moosa was partially supported by an NSERC postdoctoral fellowship}\thanks{Thomas Scanlon was partially supported by NSF Grant DMS-0071890}

\date{16 May 2002}

\begin{abstract}
Motivated by the problem of determining the structure of integral
points on subvarieties of semiabelian varieties defined over
finite fields, we prove a quantifier elimination result for
certain modules over finite simple extensions of the integers
given together with predicates for orbits of the distinguished
generator of the ring.
\end{abstract}

\maketitle

\section{introduction}
\label{intro}

The Mordell-Lang conjecture asserts that for $G$ a
semiabelian variety over the complex numbers, $X \subset G$
a subvariety, and $\Gamma \leq G(\CC)$ a finitely
generated subgroup of the complex points, the set of
points $X(\CC) \cap \Gamma$ is a finite union of
cosets of subgroups of $\Gamma$.
This fails when $\CC$ is
replaced by a field of positive characteristic.
For example, suppose that
$X$ and $G$ are defined over a finite field $\FF_q$ and
that $F: G \to G$ is the corresponding Frobenius morphism.
Let
$K:=\FF_q(X)$ and let $\Gamma \leq G(K)$ be the
$\ZZ[F]$-submodule generated by $\gamma:=\id_X:X \to X$
thought of as an element of $X(K)$.
Then $X(K) \cap \Gamma$
contains the infinite set $\{ F^n\gamma : n \in\NN \}$.
If $X$ contains
no translates of algebraic subgroups of $G$, then this fact already
contradicts the na\"{\i}ve translation of the Mordell-Lang conjecture
to positive characteristic.

Hrushovski salvages the Mordell-Lang conjecture in positive characteristic
by proving a function field version in which varieties defined over
finite fields are treated as exceptions to the general rule~\cite{hrushovski}.
In this
paper we generalize the rule so that these varieties are no longer
exceptional.

In the case that $X$ is a curve, these Frobenius orbits are the only obstruction to a clean statement of
Mordell's conjecture.  Samuel showed that if $C$ is a curve of geometric
genus at least two defined over
a finite field, then for any finitely generated field $K$ extending
the field of definition of $C$, the set of $K$-rational points on $C$ is a
finite set of Frobenius orbits~\cite{Sam}.  Continuing
with the example at the end of the first paragraph, suppose that
$Y = \overline{X + X}$ also contains no translates of infinite algebraic
subgroups of $G$.
This is the case, for example, when $X$ is a curve of
genus
at least three embedded into its Jacobian $G = J_X$.  Then $Y(K) \cap
\Gamma$
contains the set $\{F^m \gamma + F^n \gamma : n,m \in {\mathbb N} \}$.
We show in Section \ref{intended-application} of this paper that to handle the general case of the Mordell-Lang problem for
$G$ defined over a finite field we need
only permit such sums of finitely many Frobenius orbits
(together with groups) into the description of $X(K) \cap \Gamma$.
It seems that the main step in proving this result was to recognize the
correct form of these intersections.

The Mordell-Lang conjecture and related statements have model theoretic
interpretations.  On the face of it, the Mordell-Lang conjecture in its
original form may be rephrased as \emph{The structure induced on
a finitely generated subgroup of a semiabelian variety of the
complex numbers from the field structure is weakly normal.} Moreover,
as Pillay showed, the finitely generated group is actually stably embedded.
More precisely, if $K$ is an algebraically closed field of characteristic
zero and $\Gamma$ is a finitely generated subgroup of the $K$-points of
some semiabelian variety over $K$, then the theory of the structure
$(K,+,\times,\Gamma)$ is stable and the formula ``$x\in\Gamma$'' is weaky normal~\cite{pillay}.
As observed by the second author in~\cite{Sca-uniform} this result implies
a version of uniformity for the Mordell-Lang conjecture.

The present paper addresses the question of determining the model theoretic
properties of the structure induced on $\Gamma$ by $(K,+,\times,\Gamma)$,
when $K$ is an algebraically closed field of positive characteristic and
$\Gamma$ is a finitely generated Frobenius submodule of a semiabelian
variety defined over a finite field.  As certain infinite Frobenius orbits
may be definable in $\Gamma$, the induced structure cannot be
weakly normal.  However, we show that it is stable, and hence, as in Pillay~\cite{pillay}, $(K,+,\times,\Gamma)$ is stable.  As a consequence
of this analysis we obtain a uniform version of the Mordell-Lang conjecture
for semiabelian varieties over a finite field.

While the structure of integral points on semiabelian varieties defined
over
finite fields serves as motivation, we perform our technical work in the
abstract setting of $F$-spaces.  We work with a fixed finite simple
extension
$R$ of ${\mathbb Z}$, which we denote by $\ZZ[F]$,
and a class of finitely generated $R$-modules which we call $F$-spaces
(see Section \ref{set-up} for the formal definitions).  For an $F$-space $M$, an $F$-subspace $N\leq M$, a finite tuple
$a_1, \ldots, a_n \in M$, and positive integers $\delta_1, \ldots, \delta_n$;
we associate the set $S(\a, \overline{\delta}, H) := \{ \sum_{i=1}^n
F^{m_i \delta_i} a_i \} + H$.
We call such sets $F$-sets, and prove quantifier elimination and stability for
$F$-spaces with predicates for these $F$-sets.
The key is a translation
between
properties of $F$-sets and of sets definable in the trivial,
strongly minimal  structure $({\mathbb N}, \sigma, 0)$
where $\sigma$ is the successor operation.

\smallskip

Our collaboration on this paper began during the MSRI Model Theory of
Fields
program during the Spring of 1998.  Both authors thank MSRI for providing
excellent working conditions.

\bigskip
\bigskip
\section{Set-Up and Statement of Results}
\label{set-up}

Fix a unital ring $R$ that is generated by a single distinguished element $F\in R$, over the integers $\ZZ$.  We write $R=\ZZ[F]$.

\smallskip

\begin{definition}
\label{f-space}
An {\em $F$-space} is a finitely generated $R$-module, $M$, such that
$$ F^\infty M= \bigcap_{n=0}^\infty F^nM$$
is finite.
By an {\em $F$-subspace} of $M$, we mean an $R$-submodule $H\subset M$, such that the quotient module $M/H$ is again an $F$-space.
\end{definition}

\smallskip

\begin{remark}
For fixed $R=\ZZ[F]$, the following observations are immediate consequences of the definitions:
\begin{itemize}
\item The class of $F$-spaces is closed under taking products and passing to quotients by $F$-subspaces.
\item The class of $F$-subspaces is closed under taking intersections, and products.
\item The diagonal $\Delta\subset M^2$, and the graph of addition $\Sigma\subset M^3$, are $F$-subspaces.
\item If $H\subset M$ is an $F$-subspace and $\pi\colon M\to M/H$ is the quotient map, then the class of $F$-subspaces (of $M$ and $M/H$) are closed under $\pi$ and $\pi^{-1}$.
\end{itemize}
\end{remark}

\smallskip

\begin{example}
\label{intended}
We describe our intended example.
Fix $G$ a semiabelian variety over a finite field ${\mathbb F}_q$ of characteristic $p>0$.
The group variety $G$ admits an algebraic endomorphism $F:G \to G$ induced 
from the $q$-power Frobenius.
Let $R=\ZZ[F]$ be the subring of the endomorphism ring of $G$ generated by $F$.
Now, let $K$ be a finitely generated regular extension field of ${\mathbb F}_q$.  
If $\Gamma \subset G(K)$ is a finitely generated
$R$-submodule of $G(K)$, then as
$$F^{\infty}G(K)=\bigcap_{n=0}^\infty F^nG(K)=\bigcap_{n=0}^\infty G(K^{q^n})=G(\FF_q),$$
we see that $\Gamma$ is an $F$-space.

A particularly relevant case is when $\Gamma$ is the set of rational or integral points on $G$.
More precisely, if $G$ is an abelian variety, then $G(K)$ is itself a finitely generated group, and hence a finitely generated $R$-module.
We could take $\Gamma$ to be $G(K)$.
More generally, for semiabelian varieties, we could take $\Gamma$ to be $G(\R)$, where $\R\subset K$ is a finitely generated ring extension of $\FF_q$.
\end{example}

\smallskip

The rather forced notion of an $F$-subspace in Definition \ref{f-space} requires some explanation.
The problem is that a quotient of an
$F$-space by an arbitrary submodule need not be an $F$-space -- we may lose control over the infinitely $F$-divisible points.
The following observation, however, gives us a class of examples (including our intended example \ref{intended} above) where every finitely generated $R$-module is an $F$-space.
It is not hard to see that in this case, every $R$-submodule of an $F$-space will be an $F$-subspace.

\smallskip

\begin{proposition}
\label{module=space}
Suppose $R=\ZZ[F]$ is a finite extension of $\ZZ$, $F$ is not a zero-divisor, and $\displaystyle \bigcap_{n\geq 0}F^nR=0$.
Then every finitely generated $R$-module is an $F$-space.
\end{proposition}

\proof
Let $M$ be a finitely generated $R$-module, and let $N=F^\infty M$.
We wish to show that $N$ is finite.
Let $(F)$ be the ideal in $R$ generated by $F$, and let $\pi\colon R\to R/(F)$ be the quotient map.
As $N$ is finitely generated, and $FN=N$, Nakayama's Lemma implies that there is an $r\in R$ with $rN=0$ and $\pi(r)= 1$.

We claim that $r$ is not a zero-divisor.
Indeed, suppose $sr=0$ for some $s\neq 0$ in $R$.
Then $0=\pi(sr)=\pi(s)\pi(r)=\pi(s)$.
Hence, $s\in (F)$.
On the other hand, as $F^\infty R=0$, $s=F^nt$ for some $n\geq 0$ and $t\in R\setminus (F)$.
Hence, $F^ntr=0$, and as $F$ is not a zero divisor, $tr=0$.
But then, as before, $\pi(t)=0$, implying that $t\in(F)$ -- which is a contradiction.

It follows that $I=(r)\cap\ZZ$ is a nonzero ideal of $\ZZ$.
Indeed, as $R$ is a finite extension of $\ZZ$, $r$ is integral over $\ZZ$.
The minimal polynomial of $r$ over $\ZZ$ must have a nonzero constant term, since $r$ is not a zero-divisor.
This constant term is visibly in $(r)\cap\ZZ$.
Hence $\ZZ/I$ is a finite ring, and $R/(r)$, being a finitely generated $\ZZ/I$-module, is also finite.
Finally, as $r$ kills $N$, $N$ is a finitely generated $R/(r)$-module -- and hence is itself finite.
\qed

\begin{remark}
\label{intended-good}
If $G$ and $R=\ZZ[F]$ are as in Example \ref{intended}, then $R$ and $F\in R$ satisfy the conditions of Proposition \ref{module=space}.
Indeed, the endomorphism ring of $G$, and hence $R$, is a finite extension of $\ZZ$.
Since $F$ is injective on $G$ it is not a zero-divisor.
Moreover, the only infinitely $F$-divisible element of $R$ is the zero map.
To see this, choose a finitely generated field extending $\FF_q$, $L$, such that $G(L)$ is Zariski dense in $G$.
Note that every endomorphism in $R$ is defined over $\FF_q$, and hence over $L$.
If $\alpha\in F^\infty R$, then $\displaystyle \alpha G(L)\subset\bigcap_{n>0}F^nG(L)=G(k)$, where $\displaystyle k=\bigcap_{n>0}L^{q^n}$ is a finite field.
Hence $\alpha$ takes a Zariski dense subgroup of $G$ to a finite group, which implies that $\alpha$ must be the zero map.
\end{remark}

\medskip

Fix $R=\ZZ[F]$ and an $F$-space $M$.
We are interested in subsets of $M$ that are obtained as sums of orbits of points under powers of $F$.
In order to include $0$ in any such orbit, we let $\NN_*=\NN\cup\{-\infty\}$, and use the convention that $F^{-\infty}a=0$ for all $a\in M$.

\smallskip

\begin{definition}
\label{f-sets}
Suppose $M$ is an $F$-space, $\a=(a_1,\dots ,a_n)$ is a tuple from $M$, and $\overline{\delta}=(\delta_1,\dots ,\delta_n)$ is a tuple from $\NN$ with each $\delta_i>0$.
\begin{itemize}
\item[1.] A set of the form
$\displaystyle S(\a;\overline{\delta})=\{\sum_{i=1}^nF^{\delta_im_i}a_i\colon (m_1,\dots ,m_n)\in\NN_*^n\}$
is called a {\em basic groupless $F$-set}.
\item[2.] A {\em basic $F$-set} is a set of the form $S(\a;\overline{\delta})+H$ where $S(\a;\overline{\delta})$ is a basic groupless $F$-set in $M$ and $H$ is an $F$-subspace of $M$.
\end{itemize}
A {\em groupless $F$-set} is a finite union of translates of basic groupless $F$-sets.
An {\em $F$-set} is a finite union of translates of basic $F$-sets.
\end{definition}

\smallskip

\begin{remark}
\label{fixed-power}
If $a\in M$ and $\delta\in\NN$, then $S(a;\delta)$ is just the $F^\delta$-orbit of $a$ in $M$.
In the case that $\delta_1=\cdots =\delta_n=\delta$, we use the abbreviation $S(\a;\delta)$ instead of $S(\a;\overline{\delta})$.
Notice that in general, if $S(\a;(\delta_1,\dots ,\delta_n))$ is a basic groupless $F$-set, then by taking $\delta$ to be the least common multiple of the $\delta_i's$, we get that $S(\a;(\delta_1,\dots ,\delta_n))$ is a finite union of basic groupless $F$-sets of the form $S(\a';\delta)$.
That is, up to finite unions, we can take all the $\delta_i's$ to be the same.
\end{remark}

\smallskip

\begin{definition}
\label{f-structure}
An {\em $F$-structure} is a pair $(M,\S)$ where $M$ is an $F$-space and $\S=\bigcup_{n\geq 0}\S_n$, where $\S_n$ is the collection of all $F$-sets in $M^n$.
We view $(M,\S)$ as a first order structure in the language of $F$-sets -- that is, where there is a predicate for each member of $\S$.
\end{definition}

\smallskip

Here are the main results of this paper:

\theoremstyle{plain}
\newtheorem*{theorem-a}{Theorem A}
\begin{theorem-a}
The theory of an $F$-structure admits quantifier elimination and is stable {\em (Theorems \ref{qe} and \ref{stable}, respectively)}.
\end{theorem-a}

Based on an argument of the second author, we obtain the following version of Mordell-Lang for semiabelian varieties over finite fields (this is Theorem~\ref{ml}).
Fix $\Gamma\subset G(K)$ as in Example~\ref{intended}.

\theoremstyle{plain}
\newtheorem*{theorem-b}{Theorem B}
\begin{theorem-b}
If $X \subseteq G$ is a closed subvariety, then $X\cap\Gamma$ is an $F$-set.
The basic $F$-sets that appear have the form $S(a_1, \ldots, a_n;1) + (H\cap\Gamma)$ for some algebraic subgroup $H \leq G$ over $\FF_q$ and points $a_1, \ldots, a_n \in \Gamma$.
\end{theorem-b}

Combining these, we are able to conclude (this is Corollary~\ref{stable-predicate}):

\theoremstyle{plain}
\newtheorem*{theorem-c}{Theorem C}
\begin{theorem-c}
If $\U$ is an algebraically closed field extending $K$, then the theory of $(\U,+,\times,\Gamma)$ is stable.
\end{theorem-c}

As an application of our analysis, we also obtain the following uniform version version of Theorem $B$ (this is Corollary~\ref{uniform-ml}):

\theoremstyle{plain}
\newtheorem*{theorem-d}{Theorem D}
\begin{theorem-d}
Suppose $\{ X_b \}_{b \in B}$ is an algebraic family of closed subvarieties of $G$.
Then there are basic groupless $F$-sets $T_1, \ldots, T_m \subseteq \Gamma$ and algebraic subgroups $H_1, \ldots, H_m \leq G$ over $\FF_q$, such that for any $b \in B$ there is an $I \subseteq \{1, \ldots, m \}$ and points $(\gamma_i)_{i\in I}$ from $\Gamma$, such that
$$X_b\cap \Gamma = \bigcup_{i \in I} \gamma_i + T_i + (H_i \cap \Gamma).$$
\end{theorem-d}

\smallskip

One final comment.
A careful inspection of the proof of Theorem $A$ would yield quantifier elimination and stability in certain restricted languages.
For example, suppose $\displaystyle \A=\bigcup_{n\geq 0}\A_n$, where $\A_n$ is a collection of $F$-subspaces of $M^n$, such that $(M,\A)$ is an abelian structure admitting quantifier elimination.
That is, $\A$ contains $0$, $M$, the diagonal $\Delta\subset M^2$, and the graph of addition $\Sigma\subset M^3$; and $\A$ is closed under taking products, intersection, and coordinate projections.
For each $n\geq 0$, denote by $\S\A_n$ the set of all $F$-sets in $M^n$ where the $F$-subspaces that appear come from $\A_n$; and let $\displaystyle \S\A=\bigcup_{n\geq 0}\S\A_n$.
We say that $(M,\S\A)$ is the {\em $F$-structure on $(M,\A)$}.
It will follow from our proofs of Theorems \ref{qe} and \ref{stable}, that $\th(M,\S\A)$ admits quantifier elimination and is stable.

\smallskip

\begin{example}
\label{intended-f-structure}
Let $\Gamma\leq G(K)$ be as in Example \ref{intended}.
We may regard $\Gamma$ as an abelian structure by taking as basic definable subgroups of $\Gamma^n$  those groups of the form $H \cap \Gamma^n$ for $H \leq G^n$ an algebraic subgroup.
As $H \cap \Gamma^n = (\overline{H \cap \Gamma^n}) \cap \Gamma^n$, we need only consider those algebraic subgroups which are defined over $K$.
As every algebraic subgroup of $G^n$ is defined over ${\mathbb F}_q^{\alg}$ and as $K\cap\FF_q^{\alg}=\FF_q$, we need only consider those algebraic subgroups which are defined over $\FF_q$.
As such, we see that these basic definable groups are $\ZZ[F]$-submodules of $\Gamma^n$.
By Proposition \ref{module=space} and Remark \ref{intended-good}, they are in fact $F$-subspaces of $\Gamma^n$.
We obtain $\A_n$ by closing off under coordinate projections $\Gamma^{m+n}\to \Gamma^n$.
Thus $\A_n$ is a collection of $F$-subspaces of $\Gamma^n$, and letting $\displaystyle \A=\bigcup_{n\geq 0}\A_n$, we have that $(\Gamma,\A)$ is an abelian structure that admits quantifier elimination.
The $F$-structure $(\Gamma,\S\A)$ on $(\Gamma,\A)$, will then also admit quantifier elimination and be stable.
\end{example}

\bigskip
\bigskip

\section{Definability in $\NN_*$ and Groupless $F$-sets}
\label{exponents}

Recall that $\NN_*=\NN\cup\{-\infty\}$.
We also declare: {\em for all $m\in\NN_*$ and $\delta\in\NN$, $m(-\infty)=m+(-\infty)=-\infty$; $-\infty<\delta$, and $-\infty\equiv m\mod\delta$.}

We view $\NN_*$ as a structure in the language $\L^*$ where there are constant symbols for $0$ and $-\infty$, a function symbol $\sigma$ for the successor function, and unary predicates $P_{\delta}(x)$ that signify ``$x\equiv 0\mod\delta$'', for each $0<\delta\in\NN$.
When working with groupless $F$-sets, one is lead to consider certain distinguished classes of definable sets in the structure $(\NN_*,\sigma, 0,-\infty,(P_{\delta})_{\delta>0})$.

\begin{definition}
\label{basic-exponent}
By an {\em equation in $\L^*$} we mean a formula of the form $X=\sigma^r(Y)$, $X=p$, or $X=q\mod\delta$; where $X$ and $Y$ are (singleton) variables, $r\in\NN$, $p,q\in\NN_*$, and $\delta\in\NN$ is positive.
\begin{itemize}
\item 
A {\em variety} in $\NN_*^n$, is a solution set to finitely many equation in the variables $X_1,\dots X_n$.
\item
A translate of a variety is called a {\em basic set}.
\item
A {\em closed} set is a finite union of basic sets.
\end{itemize}
\end{definition}

\begin{remark}
\label{basic-operations}
Note that while varieties are closed under taking intersections, they are not closed under projections.
The class of basic sets, on the other hand, is closed under intersections, translations, and projections.
The closed sets can also be described as the projections of positive quantifier-free definable sets.
\end{remark}

\smallskip

It is convenient to work with a single congruence predicate at a time.
Given $\delta>0$, we let $\L^\delta\subset\L^*$ be the finite sublanguage consisting of $P_{\delta}(x)$ together with $\sigma$, $0$, and $-\infty$.
We say that $B\subset\NN_*^n$ is  a {\em $\delta$-variety} if it is a variety given by equations only involving congruence modulo $\delta$.
Similarly for $\delta$-basic and $\delta$-closed sets.
We gather some facts about definable sets in $(\NN_*,\sigma,0,-\infty,P_{\delta})$, for future use.

\begin{lemma}
\label{exponent-qe}
Fix a positive integer $\delta$.
\begin{itemize}
\item[(a)]
The structure $(\NN_*,\sigma,0,-\infty,P_{\delta})$ is of Morley rank one.
\item[(b)]
Every $\L^{\delta}$-definable set is a finite boolean combination of $\delta$-varieties.
\item[(c)] If $B$ and $C$ are $\delta$-varieties, with $C$ a proper subset of $B$, then either $\rm(C)<\rm(B)$ or $\dm(C)<\dm(B)$.
\item[(d)] Suppose $\{X_{\r}\colon\r\in\NN_*^m\}$ is a uniformly $\L^{\delta}$-definable family of definable subsets of $\NN_*^n$.
There is a finite union of $\delta$-varieties, $B\subset\NN_*^{m+n}$, such that for all $\r\in\NN_*^m$ for which $X_{\r}$ is nonempty, $X_{\r}\subset B_{\r}$, $\rm(X_{\r})=\rm(B_{\r})$, and $\dm(X_{\r})=\dm(B_{\r})$.
\end{itemize}
\end{lemma}

\proof
First of all, we deal with the case of $\delta=1$.
Part $(a)$ follows from the fact that $(\NN,\sigma,0)$ is strongly minimal;
 and part $(b)$ from the fact that $(\NN,\sigma,0)$ admits quantifier elimination even when the function $\sigma$ is replaced by its graph.
For part $(c)$, notice that every variety in $(\NN_*,\sigma,0,-\infty)$ is of Morley degree $1$, and that a proper subvariety is of strictly smaller Morley rank.

Consider part $(d)$ (still in the case of $\delta=1)$).
Suppose $\phi(x,y)$ is a formula such that for all $\r\in\NN_*^m$, $\phi(\r,y)$ defines $X_{\r}$.
Using part $(b)$, $\phi(x,y)$ is a boolean combination of equations.
Putting this in normal form, and then working with each disjunct at a time, we may assume that $\phi(x,y)$ is a conjunction of equations and inequations.
That is, $\phi(x,y)$ defines a set of the form $B\setminus C$, where $B\subset\NN_*^{m+n}$ is a variety and $C\subset Y$ is a finite union of varieties.
For each $\r\in\NN_*^m$, $X_{\r}=B_{\r}\setminus C_{\r}$.
Visibly, $B_{\r}$ is again a variety, and so, as we pointed out above, $C_{\r}$ is of strictly smaller Morley rank than that of $B_{\r}$ (assuming that $X_{\r}$ is nonempty).
Hence, $\rm(X_{\r})=\rm(B_{\r})$, and they are both of degree $1$.
This completes the proof of the lemma in the case of $\delta=1$.

Now suppose that $\delta>0$ is arbitrary.
For each $0\leq i<\delta$, let $\Delta_i\subset\NN_*$ be the set of points that are equivalent to $i\mod\delta$.
Then $\NN_*=\bigcup_i\Delta_i$, and by our conventions, $\Delta_i\cap\Delta_j=\{-\infty\}$ for $i\neq j$.
Moreover, each $\Delta_i$ is a $\delta$-variety.
The induced $\L^{\delta}$-structure on $\Delta_i$ is exactly $(\Delta_i,\sigma^{\delta},i,-\infty)$ -- which is a definable copy of $(\NN_*,\sigma,0,-\infty)$.
It is now not hard to see that the lemma follows from the case of $(\NN_*,\sigma,0,-\infty)$ dealt with above.
\qed

\bigskip

{\em For the rest of this section, we fix an $F$-space $M$.}

\begin{definition}
Fix $\a=(a_1,\dots ,a_n)\in M^n$.  We will use the following notation:
\begin{itemize}
\item[1.]
For $B\subset\NN_*^n$,
$$F^B\a=\{\sum_{i=1}^nF^{b_i}a_i\colon (b_1,\dots ,b_n)\in B\}\subset M.$$
\item[2.]
For $\r=(r_1,\dots ,r_n)\in\NN_*^n$,
$$F^{\r}\a=\sum_{i=1}^nF^{r_i}a_i\in M.$$
\item[3.]
For $\r,\s\in\NN_*^n$, $\r$ is {\em $\a$-equivalent} to $\s$, written $\r\sim_{\a}\s$, if $F^{\r}\a=F^{\s}\a$.
\item[4.]
For $b\in M$,
$$\log_{\a}b=\{\r\colon b=F^{\r}\a\}\subset\NN_*^n.$$
\end{itemize}
\end{definition}

\begin{remark}
The set $\log_{\a}b$ describes the ways in which $b$ can be written as sums of iterates of $F$ applied to $a_1,\dots ,a_n$.
Note that $\sim_{\a}$ is an equivalence relation and that $B\subset\NN_*^n$ is an $\a$-equivalence class if and only if $B=\log_{\a}b$ for some $b\in M$.
\end{remark}

\smallskip

It is not hard to see that if $B$ is closed, then $F^B\a$ is a groupless $F$-set.
For example, suppose $B$ is a variety of one of the following forms:
\begin{itemize}
\item[I] $\NN_*^t\times\{p\}\times\NN_*^{n-t-1}$, for some $p\in\NN_*$.
\item[II] $\{(m_1,\dots ,m_n)\in\NN_*^n\colon m_s=\sigma^r(m_t)\}$, for some $r\in\NN$.
\item[III] $\{(m_1,\dots ,m_n)\in\NN_*^n\colon m_t\equiv q\mod \delta\}$, for some $q\in\NN_*$ and $\delta\in\NN$ positive.
\end{itemize}
Then $F^B\a$ is a translate of a basic groupless $F$-set of the form (respectively):
\begin{itemize}
\item[I] $F^pa_t+S((a_1,\dots ,a_{t-1},a_{t+1},\dots ,a_n);1)$
\item[II] $S((a_1,\dots ,a_{s-1},F^ra_{s}+a_{t},a_{s+1},\dots ,a_{t-1},a_{t+1},\dots ,a_n);1)$
\item[III] $S((a_1,\dots ,a_{t-1},F^qa_t,a_{t+1},\dots ,a_n);(1,\dots,\delta ,\dots ,1))$
\end{itemize}
Moreover, if $B_2=\overline{r}+B_1$ then
$$F^{B_2}(a_1,\dots ,a_n)=F^{B_1}(F^{r_1}a_1,\dots ,F^{r_n}a_n).$$
Finally, if $\displaystyle B=\bigcup_{i=1}^lB_i\subset\NN_*^n$, then $\displaystyle F^B\a=\bigcup_{i=1}^lF^{B_i}\a$.
In fact we have:

\begin{lemma}
\label{closed-exponent=groupless}
A subset $S\subset M$ is a groupless $F$-set if and only if it is of the form $F^B\a$ for some tuple $\a\in M^n$ and closed $B\subset\NN_*^n$.
\end{lemma}

\proof
It is clear from the above discussion, that if $B$ is closed then $F^B\a$ is a groupless $F$-set.  Now suppose $S\subset M$ is a groupless $F$-set.  If it is a basic groupless $F$-set of the form $S(\a;\overline{\delta})$, then $S=F^B\a$, where $B$ is the variety set $\{(r_1,\dots ,r_n)\colon r_i\equiv 0\mod\delta_i\}$.
If $S=c+F^B\a$, with $B$ closed, then we can write $S$ as $F^{(0\times B)}(c,\a)$.
Finally, suppose $S=F^{B_1}\a\cup F^{B_2}\b$.
Let
$$B_1^\prime=B_1\times\{-\infty\}\times\cdots\times\{-\infty\}$$
and
$$B_2^\prime=\{-\infty\}\times\cdots\times\{-\infty\}\times B_2,$$
where the number of $-\infty$'s that we have attached corresponds to the arity of $\a$ and $\b$ respectively.
Let $\c$ be the concatenation of $\a$ and $\b$.
Then $S=F^{(B_1^\prime\cup B_2^\prime)}\c$, and we see that every groupless $F$-set is of this form.
\qed

\medskip

As in the case of closed sets, there is a relationship between arbitrary definable sets in $\NN_*$ and finite boolean combination of groupless $F$-sets in $M$.
However, in order to see this correspondence, we need to consider $\a$-equivalence.

\begin{proposition}
\label{log-closed}
Suppose $\a\in M^n$.
For all $b\in M$, $\log_{\a}b\subset\NN_*^n$ is closed.
\end{proposition}

The proof of Proposition \ref{log-closed} is somewhat technical and we delay it for the time being.
In Theorem \ref{log} below, we will not only prove that $\log_{\a}b$ is closed, but we will also describe how the set varies with $b$.

\begin{corollary}
\label{equivalence}
Suppose $\a\in M^n$.
Then $\a$-equivalence is a {\em (}closed{\em )} definable equivalence relation in $\NN_*$.
\end{corollary}

\proof
Notice that $\r\sim_{\a}\s$ if and only if $F^{\r}\a-F^{\s}\a=0$.
Letting $\a^\prime\in M^{2n}$ be the tuple $(\a,-\a)$, this is equivalent to
$(\r,\s)\in\log_{\a^\prime}0$.
By Proposition \ref{log-closed}, $\log_{\a^\prime}0$ is a closed set.
\qed

\smallskip

\begin{proposition}
\label{definable-exponent=qf}
Suppose $M$ is an $F$-space, $\a\in M^n$, and $Y\subset\NN_*^n$ is definable.
Then $F^Y\a$ is a finite boolean combination of groupless $F$-sets.
Conversely, if $U$ is a finite boolean combination of groupless $F$-sets, then for some tuple $\a\in M^n$, and definable $Y\subset\NN_*^n$, $U=F^Y\a$.
\end{proposition}

\proof
Suppose $Y\subset\NN_*^n$ is definable and let $\delta>0$ be such that $Y$ and $\sim_{\a}$ are $\L^\delta$-definable.
Working in the language $\L^\delta$, we proceed by induction on the Morley rank of $Y$.
The case of $\rm Y=0$ is trivial.
Let $r=\rm Y>0$.
Taking finite unions we may assume that the Morley degree of $Y$ is $1$.
By Lemma \ref{exponent-qe}, $Y$ is contained in a $\delta$-closed set $B$ of Morley rank $r$ and Morley degree $1$.
(In fact, this is a special case of part $(d)$ of that lemma; we are considering a single definable set rather than a family of them.)
It follows that $Y=B\setminus C$, where $C\subset B$ is an $\L^{\delta}$-definable set of Morley rank strictly less than $r$.

Let $D\subset C$ be the set of all $\r\in C$ such that if $\s\in B$ with $\s\sim_{\a}\r$, then $\s\in C$.
Then $F^Y\a=F^B\a\setminus F^D\a$.
Note that $D$ is $\L^\delta$-definable, and as it is contained in $C$, it is also of Morley rank strictly less than $r$.
By induction, $F^D\a$ is a boolean combination of groupless $F$-sets.
As $B$ is $\delta$-closed, we already know that $F^B\a$ is a groupless $F$-set.
Hence, $F^Y\a$ is a boolean combination of groupless $F$-sets.

For the converse, it is suffcient to show that if $\a\in M^n$ and $C\subset B$ are closed sets in $\NN_*^n$, then $F^B\a\setminus F^C\a=F^Y\a$ for some definable set $Y\subset\NN_*^n$.
Indeed, $Y$ will be the set of all $\r\in B$ such that for no $\s\in C$ is $\r\sim_{\a}\s$.
This is visibly definable.
\qed

\smallskip

Suppose $G\subset M$ is an $F$-subspace and $\pi\colon M\to M/G$ is the quotient map.
From the definition of $F$-sets it is clear that the image or preimage of an $F$-set under $\pi$ is again an $F$-set.
Moreover, as $\pi^{-1}$ respects boolean operations, finite boolean combinations of $F$-sets are also preserved under $\pi^{-1}$.
To show that such boolean combinations are preserved under $\pi$ amounts to quantifier elimination for $(M,\S)$, and will take some work.
However, for groupless $F$-sets, the situation is now much easier:

\begin{corollary}
\label{groupless-qe}
Suppose $G\subset M$ is an $F$-subspace and $\pi\colon M\to M/G$ is the quotient map.
The image under $\pi$ of a boolean combination of groupless $F$-sets is again a boolean combination of groupless $F$-sets.
\end{corollary}

\proof
By Proposition \ref{definable-exponent=qf}, $U=F^Y\a$ where $\a\in M^n$ and $Y\subset\NN_*^n$ is definable.
But then, $\pi U=F^Y(\pi\a)$, which, by the proposition again, is a boolean combination of groupless $F$-sets.
\qed

\medskip

We introduce a usefull subclass of the quantifier-free definable sets that behave very much like $F$-sets but give us greater flexibility.

\begin{definition}
Suppose $M$ is an $F$-space.
A {\em generalised $F$-set} is a finite union of sets of the form $V+G$ where $V$ is a boolean combination of groupless $F$-sets and $G$ is an $F$-subspace of $M$.
\end{definition}

\begin{remark}
\label{generalised-qf}
Generalised $F$-sets are quantifier-free definable -- that is, they are boolean combinations of $F$-sets.
Indeed, if $V+G\subset M$ is as above, and $\pi$ is reduction modulo $G$, then $V+G=\pi^{-1}(\pi V)$.
Now by Corollary \ref{groupless-qe}, $\pi V$ is a boolean combination of groupless $F$-sets.
Hence $V+G$ is a boolean combination of $F$-sets.
\end{remark}

\begin{proposition}
\label{positive-qf}
The class of generalised $F$-sets is preserved under taking unions, intersections, and images/preimages of quotient maps.
\end{proposition}

\proof
For unions, this is part of the definition.
Let us deal with images and preimages of projections.
Suppose $\pi\colon M\to M/G$ is the quotient map where $G$ is an $F$-subspace of $M$.
Since $\pi$ and $\pi^{-1}$ commute with unions,
and as the class of $F$-subspaces is closed under $\pi$ and $\pi^{-1}$,
we need only consider boolean combinations of groupless $F$-sets.
But by Corollary \ref{groupless-qe}, if $V\subset M$ is such,  then $\pi V$ is again a boolean combination of groupless $F$-sets.
On the other hand if $V\subset M/G$ is a boolean combination of groupless $F$-sets, then write $V=F^Y\a$, with $\a$ a tuple from $M/G$ and $Y$ a definable set from $\NN_*$.
Then $\pi^{-1}(V)=F^Y\b+G$, for any $\b$ from $M$ whose reduction modulo $G$ is $\a$.
This is a generalised $F$-set.

Now for intersections.
Given tuples $\a,\b$ from $M$, definable sets $X,Y$ from $\NN_*$, and $F$-subspaces $H, G$, we need to show that $(F^X\a+G)\cap(F^Y\b+H)$ is a generalised $F$-set.
We deal with cases seperately:

{\em Case $1$: $G=0$ and $F^Y\b=0$}.
In this case we have that $F^X\a\cap H= F^Z\a$ where $Z=\log_{\pi\a}0\cap X$, where $\pi$ denotes reduction modulo $H$.
As $Z$ is definable, $F^X\a\cap H$ is a boolean combination of groupless $F$-sets.

{\em Case $2$: $F^Y\b=0$}.
Letting $\pi$ be the reduction modulo $G\cap H$ map, we see that
$(F^X\a+G)\cap H=\pi^{-1}[(\pi(F^X\a)+\pi G)\cap \pi H]$.  Hence we may assume that $G\cap H=0$.
Moreover, $(F^X\a+G)\cap H=[(F^X\a\cap(G+H))+G]\cap H$.
By case $1$, $F^X\a\cap(G+H)$ is a boolean combination of groupless $F$-sets.
Hence, we may work entirely inside $G+H$ -- that is, we may assume that $M=G+H$.
But by these reductions, $G+H=G\oplus H$, so we have a projection map $\eta\colon G+H\to H$.
Moreover, $(F^X\a+G)\cap H=\eta(F^X\a)$, which we already know is a boolean combination of groupless $F$-sets.

{\em Case $3$: The general case.}
Note that $(F^X\a+G)\cap(F^Y\b+H)$ is the projection onto the first coordinate of $[(F^X\a+G)\times (F^Y\b+H)]\cap\Delta$ where $\Delta$ is the diagonal.
On the other hand, it is not hard to see that $(F^X\a+G)\times (F^Y\b+H)$ is again a generalised $F$-set.
Intersecting with the diagonal -- which is an $F$-subspace -- gives us another generalised $F$-set, by case $2$.
So $(F^X\a+G)\cap(F^Y\b+H)$ is the projection of a generalised $F$-set and hence is also a generalised $F$-set.
\qed

\begin{remark}
\label{positive-qf-special}
The above arguments also show that the class of $F$-sets is preserved under intersections (we already know that they are preserved under unions and images/preimages of quotient maps).
Indeed, an $F$-set is a union of sets of the form $F^B\a+G$ where $B$ is a closed set from $\NN_*$ and $G$ is an $F$-subspace of $M$.
Using the fact that the $\log$ sets are closed (and not just definable), it is not hard to see that the proof of the above proposition also works for this class.
\end{remark}

\bigskip
\bigskip

\section{Uniform Definability of Logarithmic Equivalence}
\label{technicalities}

Fix an $F$-space $M$.
In this section, we aim to give a uniform description of sets in $\NN_*$ of the form $\log_{\a}b$, which will in particular prove Proposition \ref{log-closed}.
In order to do so, it is convenient to introduce additional terminology that distinguishes certain kinds of translates of basic sets in $\NN_*$:

\begin{definition}
\label{direct-sum}
Suppose $B\subset\NN_*^n$ is closed.
A tuple $\overline{r}\in\NN^n$ is called {\em disjoint from $B$} if for all $i\leq n$, if the $i$th coordinate of $\overline{r}$ is nonzero then every element of $B$ has zero as its $i$th coordinate.
In this case, we write $\overline{r}\oplus B$ to mean the translate $\overline{r}+B$.
We call $\overline{r}\oplus B$ a {\em disjoint translate} of $B$.
\end{definition}

\smallskip

The purpose of introducing disjoint translates is that if $B\subset\NN_*^n$ is closed and $\r\oplus B$ is a disjoint translate of $B$, then we have a particularly simple description of $F^{(\r\oplus B)}\a$ in terms of $F^{B}\a$.
Let $I_B\subset\{1,\dots n\}$ be the set of indices $i$ such that some element of $B$ has a nonzero $i$th coordinate.
We associated to $B$ a {\em canonical contraction}, denoted by $B^\circ$, that is obtained from $B$ by replacing the $j$th coordinate of every element of $B$ with $-\infty$, for all $j\notin I_B$.
Let $\r^\circ$ be obtained from $\r$ by replacing the $i$th coordinate of $\r$ with $-\infty$ for all $i\in I_B$.
It is not hard to see that
$$F^{(\r\oplus B)}\a=F^{\r^\circ}\a+F^{B^\circ}\a.$$
Notice that $F^{B^\circ}\a$ does not depend on $\r$, only the point by which $F^{B^\circ}\a$ is being translated does.
That is, as you disjoint translate a closed set, the groupless $F$-sets you obtained only vary by translation.

Another advantage of disjoint translates is that they form a uniformly definable family of sets.
Suppose $B\subset\NN_*^n$ is definable.
Let $C$ be the set of $\r\in\NN_*^n$ such that $\r$ is disjoint from $B$.
Then,
$$\{\r\oplus B\colon \r\in C\}$$
is a uniformly definable family of basic sets parameterised by $C$.
Notice that since we do not have ``$+$'' in our language $\L^*$, arbitrary translations of a definable set do {\em not} form a uniformly definable family.
The following lemma says that every unformly definable family of varieties is essentially of this form.

\begin{lemma}
\label{uniform-variety=disjoint}
Suppose $V\subset\NN_*^{m+n}$ is a variety.
There is a variety $W\subset\NN_*^n$ such that for all $\r\in\NN_*^m$, $V_{\r}$ is either empty or a disjoint translate of $W$.
\end{lemma}

\proof
Let us use $X_1,\dots ,X_m, Y_1,\dots ,Y_n$ as coordinate variables for $\NN_*^{m+n}$.
Arrange the variables so that for some $1\leq\ell\leq n$ the equations defining $V$ imply a condition of the form $Y_i=\sigma^t(X_j)$ or $X_j=\sigma^t(Y_i)$ for some $j\leq m$,  if and only if $i\leq\ell$.
It follows that as $\r\in\NN_*^m$ varies the last $n-\ell$ coordinates of $V_{\r}$ remain fixed.
Indeed, we have that $V_{\r}=\s\times W^\prime$, where $\s\in\NN_*^{m+\ell}$ varies with $r$ and $W^\prime\subset\NN_*^{n-\ell}$ is a fixed variety.
Then $W=\{0\}^{m+\ell}\times W^\prime$ satisfies the conclusion of part $(d)$.
\qed

\bigskip

In any case, we aim to prove the following:

\begin{theorem}
\label{log}
Suppose $\a=(a_1,\dots ,a_n)\in M^n$.
There exist basic sets $B_1,\dots ,B_l$ in $\NN_*^n$,
such that for all $b\in M$, there is some $J\subset\{1,\dots ,l\}$, and for each $j\in J$ there is some $\r_j\in \NN_*^n$ disjoint from $B_j$, such that
$$\log_{\a}b=\bigcup_{j\in J}(\r_j\oplus B_j).$$
\end{theorem}

\medskip

We will proceed via a series of lemmas.
The main technical tool will be the following valuation on $M$.
Let $M^\circ=\bigcup_n\ker F^n$.
Note that since $M$ is finitely generated as an $R$-module, so is $M^\circ$, and hence $M^\circ=\ker F^N$ for some $N\geq 0$.
For each $i\geq 0$, let $M_i=M^\circ+F^iM$.
These are the points that are $F^i$ divisible modulo $M^\circ$.
We obtain a filtration of $M$, and define $M_\omega$ to be the intersection of this descending chain of $R$-submodules:
$$M_0=M\supset M_1\supset M_2\supset\dots\supset M_\omega=\bigcap_{n=0}^\infty M_n$$
This in turn induces a valuation on $M$, $v\colon M\to\omega+1$, given by $v(x)\geq n$ if and only if $x\in M_n$.

\smallskip

\begin{lemma}
For all $x\in M$, $v(Fx)=1+v(x)$.
\end{lemma}

\proof
This is clear if $v(x)=\omega$ (i.e., if $x\in M_\omega$).
Assume $v(x)=m<\omega$.
Let $y\in M$ and $\delta\in M^\circ$ be such that $x=F^my+\delta$.
Then $Fx=F^{m+1}y+F\delta$, and hence $v(Fx)\geq m+1$.
Now suppose that $v(Fx)>m+1$.
That is,  there exists $z\in M$ and $\gamma\in M^\circ$, such that
$Fx=F^{m+2}z+\gamma$.
Hence $F^{m+1}y+F\delta=F^{m+2}z+\gamma$, and so $F^{m+1}(y-Fz)\in M^\circ$.
Clearly, this implies that $y-Fz\in M^\circ$.
That is, $y=Fz+\alpha$ for some $\alpha\in M^\circ$.
It follows that $x=F^{m+1}z+\delta+F\alpha\in M_{m+1}$, contradicting the fact that $v(x)=m$.

\qed

\begin{lemma}
\label{value-sum}
Suppose $a_1,\dots ,a_n\in M\setminus M_\omega$. 
If $k_1,\dots ,k_n \in\NN$ satisfy $k_{j}-k_1>v(a_1)$ for all $j>1$;
then, $v(\sum_{i=1}^nF^{k_i}a_i)=k_1+v(a_1)$.
\end{lemma}

\proof
Notice that for each $j>1$,
$$v(F^{k_j}a_j)=k_j+v(a_j)>k_1+v(a_1)+v(a_j)\geq k_1+v(a_1)= v(F^{k_1}a_1).$$
The lemma now follows from the fact that the valuation of a sum is the unique minimum valuation of the summands (if such a unique minimum exists).

\qed

\smallskip

The following lemma is the main technical step toward proving Theorem \ref{log}.

\begin{lemma}
\label{equal-sums}
Suppose $a_1,\dots ,a_n\in M\setminus M_\omega$.
There exists $N>0$ such that for all $p\leq n$, $k_1<\dots<k_n$, and $l_1,\dots ,l_p$ from $\NN$ satisfying
\begin{itemize}
\item[(a)] $k_1>N$ and $k_i-k_{i-1}>N$ for all $1<i\leq n$; and,
\item[(b)] $l_i>N$, and $|l_i-l_j|>N$ for all $i\neq j$ from $\{1,\dots ,p\}$;
\end{itemize}
if
$$\sum_{i=1}^n F^{k_i}a_i=(\sum_{i=1}^pF^{l_i}a_i) \mod M_\omega,$$
then for some permutation $\sigma$ of $\{1,\dots ,p\}$,
and for each $i\leq p$, 
$$F^{l_{\sigma(i)}}a_{\sigma(i)}-F^{k_i}a_i\in M_\omega.$$
In particular, $p=n$ and $k_i+v(a_i)=l_{\sigma(i)}+v(a_{\sigma(i)})$ for all $i\leq n$.
\end{lemma}

\proof
We first explain how the ``in particular'' clause follows from the main conclusion of the lemma.
By the hypothesis of the lemma, 
$$\sum_{i=p+1}^n F^{k_i}a_i=\sum_{i=1}^p(F^{l_{\sigma(i)}}a_{\sigma(i)}-F^{k_i}a_i)\mod M_\omega.$$
The main conclusion of the lemma thus implies that $\sum_{i=p+1}^n F^{k_i}a_i\in M_\omega$.
By Lemma \ref{value-sum}, if $p<n$ then this sum has value $k_{p+1}+v(a_{p+1})$, and hence cannot be in $M_\omega$.
It must therefore be that $p=n$.
Also, $F^{l_{\sigma(i)}}a_{\sigma(i)}-F^{k_i}a_i\in M_\omega$ implies that the two summands must have the same value, that is $k_i+v(a_i)=l_{\sigma(i)}+v(a_{\sigma(i)})$.

\smallskip

We now begin to prove the lemma itself. Let $m=\max\{v(a_i)\colon i\leq n\}$ and
$$N=\max\{m,v(F^{m-v(a_i)}a_i-F^{m-v(a_j)}a_j)\colon F^{m-v(a_i)}a_i-F^{m-v(a_j)}a_j\notin M_\omega\}.$$
Suppose $k_1,\dots ,k_n,l_1,\dots ,l_p$ satisfy the hypotheses of the Lemma, and let $\sigma$ be a permutation of $\{1,\dots ,p\}$ so that
$l_{\sigma(1)}<\dots<l_{\sigma(p)}$.

\smallskip

We proceed by induction on $i\leq p$, dealing first with the case of $i=1$.
Assume $F^{l_{\sigma(1)}}a_{\sigma(1)}-F^{k_1}a_1\notin M_\omega$ and seek a contradiction.
Note that by Lemma \ref{value-sum}, our choice of $N$, and our choice of $\sigma$, we already know that $k_1+v(a_1)=l_{\sigma(1)}+v(a_{\sigma(1)})$.
Using this we compute:
\begin{eqnarray*}
F^{l_{\sigma(1)}}a_{\sigma(1)}-F^{k_1}a_1
& = & F^{l_{\sigma(1)}+v(a_{\sigma(1)})-m}(F^{m-v(a_{\sigma(1)})}a_{\sigma(1)})-F^{k_1}a_1 \\
& = & F^{k_1+v(a_1)-m}(F^{m-v(a_{\sigma(1)})}a_{\sigma(1)})-F^{k_1}a_1 \\
& = & F^{k_1+v(a_1)-m}(F^{m-v(a_{\sigma(1)})}a_{\sigma(1)}
-F^{m-v(a_1)}a_1)\\
\end{eqnarray*}
Evaluating both sides we get that
$$v(F^{l_{\sigma(1)}}a_{\sigma(1)}-F^{k_1}a_1)=k_1+v(a_1)-m+ v(F^{m-v(a_{\sigma(1)})}a_{\sigma(1)}-F^{m-v(a_1)}a_1).$$
Now, since we are assuming that $F^{l_{\sigma(1)}}a_{\sigma(1)}-F^{k_1}a_1\notin M_\omega$ this gives
$$ v(F^{l_{\sigma(1)}}a_{\sigma(1)}-F^{k_1}a_1) \leq k_1+v(a_1)-m+N< k_2+v(a_1)-m\leq k_2$$
If $p=1$, then as $\displaystyle \sum_{i=2}^n F^{k_i}a_i=(F^{l_{\sigma(1)}}a_{\sigma(1)}-F^{k_1}a_1)\mod M_\omega$, we have that
$$k_2+v(a_2)=v(F^{l_{\sigma(1)}}a_{\sigma(1)}-F^{k_1}a_1)<k_2,$$
which is a contradiction.
So we may assume that $p\geq 2$.
In this case, we have that
$$k_1+v(a_1)-m+N = l_{\sigma(1)}+v(a_{\sigma(1)})-m+N \leq l_{\sigma(1)}+N < l_{\sigma(2)},$$
and so $v(F^{l_{\sigma(1)}}a_{\sigma(1)}-F^{k_1}a_1)$ is also strictly bounded by $l_{\sigma(2)}$.
On the other hand, for some $c\in M_\omega$
\begin{eqnarray*}
k_2+v(a_2)
& = & v(\sum_{j=2}^nF^{k_j}a_j) \\
& = & v(-F^{k_1}a_1+\sum_{j=1}^nF^{k_j}a_j) \\
& = & v(c-F^{k_1}a_1+\sum_{j=1}^pF^{l_{\sigma(j)}}a_{\sigma(j)}) \\
& = & v(c+F^{l_{\sigma(1)}}a_{\sigma(1)}-F^{k_1}a_1+\sum_{j=2}^pF^{l_{\sigma(j)}}a_{\sigma(j)}) \\
\end{eqnarray*}
By Lemma \ref{value-sum}, $\displaystyle v(\sum_{j=2}^pF^{l_{\sigma(j)}}a_{\sigma(j)})=l_{\sigma(2)}+v(a_{\sigma(2)})$.
As $v(F^{l_{\sigma(1)}}a_{\sigma(1)}-F^{k_1}a_1)$ is strictly less than $l_{\sigma(2)}$, we again have the contradiction
$$k_2+v(a_2)=v(F^{l_{\sigma(1)}}a_{\sigma(1)}-F^{k_1}a_1)<k_2.$$
This proves that $F^{l_{\sigma(1)}}a_{\sigma(1)}-F^{k_1}a_1\in M_\omega$.
Moreover, if $p\geq 2$ then we have actually shown that $k_2+v(a_2)=l_{\sigma(2)}+v(a_{\sigma(2)})$.

\smallskip

Fixing $i\leq p$, our induction hypothesis is that $(F^{l_{\sigma(s)}}a_{\sigma(s)}-F^{k_s}a_s)\in M_\omega$, for all $s<i$, and that $k_i+v(a_i)=l_{\sigma(i)}+v(a_{\sigma(i)})$.
If $p=i=n$, then we are done.

Assuming $(F^{l_{\sigma(i)}}a_{\sigma(i)}-F^{k_i}a_i)\notin M_\omega$ we get, using exactly the same computations as for $i=1$, that $v(F^{l_{\sigma(i)}}a_{\sigma(i)}-F^{k_i}a_i)$ is strictly bounded by $k_{i+1}$; and if $p>i$, also by $l_{\sigma(i+1)}$.
If $p=i$, then as $\displaystyle \sum_{j=i+1}^n F^{k_j}a_j=\sum_{s=i}^p(F^{l_{\sigma(s)}}a_{\sigma(s)}-F^{k_s}a_s)\mod M_\omega$, and as
$(F^{l_{\sigma(s)}}a_{\sigma(s)}-F^{k_s}a_s)\in M_\omega$,
for all $s<i$, we have the contradiction
$$k_{i+1}+v(a_{i+1})=v(F^{l_{\sigma(i)}}a_{\sigma(i)}-F^{k_i}a_i)<k_{i+1}.$$
So we may assume that $p>i$.
Similarly to the $i=1$ case, we have for some $c\in M_\omega$,
$$k_{i+1}+v(a_{i+1})=v(c+\sum_{s=1}^i(F^{l_{\sigma(s)}}a_{\sigma(s)}-F^{k_s}a_s)+\sum_{j=i+1}^pF^{l_{\sigma(j)}}a_{\sigma(j)})$$
By induction $(F^{l_{\sigma(s)}}a_{\sigma(s)}-F^{k_s}a_s)\in M_\omega$, for all $s<i$.
Since
$$v(F^{l_{\sigma(i)}}a_{\sigma(i)}-F^{k_i}a_i)<l_{\sigma(i+1)}\leq v(\sum_{j=i+1}^pF^{l_{\sigma(j)}}a_{\sigma(j)})$$
we still have that $k_{i+1}+v(a_{i+1})=v(F^{l_{\sigma(i)}}a_{\sigma(i)}-F^{k_i}a_i)$ and the contradiction $k_{i+1}+v(a_{i+1})<k_{i+1}$ still follows.
This proves Lemma \ref{equal-sums}.

\qed

\smallskip

The following lemma is a special case of Theorem \ref{log}.

\begin{lemma}
\label{special-log}
Suppose $a_1,\dots ,a_n\in M_\omega$.
Then there exists a finite sequence of basic subsets of
$\NN_*^n$, $B_1,\dots ,B_l$, such that
for all $b\in M$, there is some $J\subset\{1,\dots ,l\}$ such that $\log_{\a}b=\bigcup_{j\in J}B_j$.
\end{lemma}

\proof
Let us first deal with the case when each $a_i\in F^\infty M$.
Note that $F$ is injective on $F^\infty M$, since $\bigcup_i\ker F^i$ is the kernel of some fixed power of $F$.
Let $\delta>0$ be such that $F^{\delta}a_i=a_i$ for all $i\leq n$ (such a $\delta$ exists since $F^\infty M$ is finite).
So for all $m\in\NN_*$, $F^ma_i=F^{(m'\delta+r)}a_i=F^ra_i$ where $r\in\NN_*$ is strictly less than $\delta$ and $m'\in\NN_*$ (if $m=-\infty$ then set $m'=0$ and $r=-\infty$).
Hence, if $\overline{m}\in\log_{\a}b$, then there exists $\overline{r}_\circ\in\log_{\a}b$ each of whose coordinates are strictly less than $\delta$, and $\overline{m}'\in\NN_*^n$, such that $\overline{m}=\overline{m}'(\delta,\dots ,\delta)+\overline{r}_\circ$.
Moreover, for all $\overline{m}'\in\NN_*^n$, $\overline{m}'(\delta,\dots ,\delta)+\overline{r}_\circ\in\log_{\a}b$.
That is,
$B(\overline{r}_\circ)=\overline{r}_\circ+\{\overline{n}\in\NN_*^n\colon\overline{n}\equiv 0\mod\delta\}$
is contained in $\log_{\a}b$.
Note that $B(\overline{r}_\circ)$ is a basic set.
Now, the collection $\{B(\overline{r})\}$ as $\overline{r}$ varies among all $n$-tuples each of whose coordinates are strictly less than $\delta$, is finite (independently of $b$) and covers all of $\NN^n$.
It folows that $\log_{\a}b$ is equal to the union of some of these $B(\overline{r})$'s;
and we have proved the lemma in this case.

Now consider the general case (where each $a_i$ is only assumed to be in $M_\omega$).
Let $N_1\in\NN$ be such that $M^\circ=\ker F^{N_1}$.
Note that for all $a\in M_\omega$, $F^{N_1}a\in F^\infty M$.
Fixing $i\leq n$ and $r<N_1$, the set of elements $\overline{m}\in\log_{\a}b$  satisfying $m_i=r$ can be described as the set of all
$$(m_1,\dots ,m_{i-1},r,m_{i+1},\dots ,m_n)\in\NN_*^n$$
such that
$$(m_1,\dots ,m_{i-1},m_{i+1},\dots ,m_n)\in\log_{\a^\prime}(b-F^ra_i),$$
where $\a^\prime=(a_1,\dots ,a_{i-1},a_{i+1},\dots ,a_n)$.
For the case $n=1$ this is just the singleton $\{r\}$ (a basic set).
Otherwise, we proceed by induction.
The induction hypothesis applied to $\a^\prime$ describes $\log_{\a^\prime}(b-F^ra_i)$ as a union of translates of basic subsets of $\NN_*^{n-1}$ chosen from a finite collection that depends only on $\a^\prime$ (and hence only on $\a$).
It should be clear how we obtain the elements of $\log_{\a}b$ whose $i$th coordinates are $r$ as translates of copies of these sets in $\NN_*^n$.
Since there are only finitely many possible values for $i$ and for $r$ -- independent of $b$ -- we are left to consider only those elements of $\log_{\a}b$ where $m_i\geq N_1$ for all $i\leq n$.
But this latter set is equal to $(N_1,\dots ,N_1)+\log_{\a^{\prime\prime}}b$, where $\a^{\prime\prime}=(F^{N_1}a_1,\dots ,F^{N_1}a_n)$.
As each $F^{N_1}a_i\in F^\infty M$, we have already proved the desired conclusion for $\log_{\a^{\prime\prime}}b$.
This completes the proof of Lemma \ref{special-log}.
\qed

\bigskip

Finally we are in a position to prove Theorem \ref{log} itself (and hence, in particular, Proposition \ref{log-closed}).
Let us restate the Theorem:
{\em
Suppose $\a=(a_1,\dots ,a_n)\in M^n$.
There exist basic sets $B_1,\dots ,B_l$ in $\NN_*^n$,
such that for all $b\in M$, there is some $J\subset\{1,\dots ,l\}$, and for each $j\in J$ there is some $\r_j\in \NN_*^n$ disjoint from $B_j$, such that
$$\log_{\a}b=\bigcup_{j\in J}(\r_j\oplus B_j).$$
}

\proof
Arrange the indices so that $a_1,\dots ,a_t\notin M_\omega$ and $a_{t+1},\dots ,a_n\in M_\omega$.
We proceed by induction on $t$.
The case $t=0$ is taken care of by Lemma \ref{special-log}.
Now suppose that $0<t\leq n$, and fix $b\in M$.
We will describe $\log_{\a}b$ making sure that the data that appears in this description are independent of $b$.

Let $N$ be the bound given by Lemma \ref{equal-sums} applied to $a_1,\dots ,a_t$.
We can divide $\log_{\a}b$ canonically into the following disjoint pieces:
\begin{itemize}
\item[$1$.] $\log_{\a}b^\sharp$: the set of those $\overline{m}\in\log_{\a}b$ for which $m_i$ and $|m_i-m_j|$ are either greater than $N$ or equal to $-\infty$, for all $i,j\leq t$.
\item[$2_{i,r}$.] $\{\m\in\log_{\a}b\colon m_i=r\}$, for fixed $i\leq t$ and $0\leq r\leq N$.
\item[$3_{i,j,r}$.] The set of tuples in $\{\m\in\log_{\a}b\colon m_j=m_i+r\}$ all of whose coordinates are greater than $N$, for fixed $i<j\leq t$ and $0\leq r\leq N$. 
\end{itemize}
We will analyse these pieces of $\log_{\a}b$ separately, and in the order presented above.

Fixing $\overline{m}\in\log_{\a}b^\sharp$, we may rearrange indices so that $m_1,\dots ,m_p$ are not $-\infty$ and $m_{p+1}=\dots =m_t=-\infty$.
We want to argue that there are only finiteley many possible choices for $(m_1,\dots m_t)$ in $\log_{\a}b^\sharp$.
Let $\overline{m}^\prime\in\log_{\a}b^\sharp$ be another element.
Notice that
$$\sum_{i=1}^tF^{m^\prime_i}a_i=(\sum_{i=1}^pF^{m_i}a_i) \mod M_\omega.$$
Indeed, $F^{-\infty}x=0$ for all $x\in M$, and $a_{t+1},\dots a_n\in M_\omega$.
Lemma \ref{equal-sums} now ensures that $m^\prime_{p+1}=\dots =m^\prime_t=-\infty$ as well, and that there is only one possible choice of $m_i+v(a_i)$ (and hence of $m_i$), up to permutations of $\{1,\dots ,p\}$.
It follows that among the tuples in $\log_{\a}b^\sharp$ there are only finitely many choices for $(m_1,\dots ,m_t)$; where the number of possible choices is bounded independently of $b$.

Fixing any such choice of $r_1,\dots ,r_t$, call the set of $\overline{m}\in\log_{\a}b^\sharp$ where $m_i=r_i$ for $i\leq t$, $\log_{\a}b^\sharp(\overline{r})$.
Then
$$\log_{\a}b^\sharp(\overline{r})=(r_1,\dots ,r_t,0,\dots ,0)\oplus(\{0\}\times\dots\times\{0\}\times\log_{\a^{\prime\prime}}(b-\sum_{i=1}^tF^{r_i}a_i))$$
where $\a^{\prime\prime}=(a_{t+1},\dots ,a_n)$.
Each of the coordinates of $\a^{\prime\prime}$ are in $M_\omega$, and so Lemma \ref{special-log} applies.
This describes $\log_{\a}b^\sharp(\overline{r})$ in terms of unions of disjoint translates of basic sets that come from a finite collection that does not depend on $b$.
Ranging over the possible choices of $(r_1,\dots ,r_t)$ -- which recall is bounded independently of $b$ -- we obtain the desired conclusions for $\log_{\a}b^\sharp$.

Fix $i\leq t$ and $0\leq r<N$, and consider case $2_{i,r}$.
The set of elements $\overline{m}\in\log_{\a}b$  satisfying $m_i=r$ can be described as the set of all
$$(m_1,\dots ,m_{i-1},r,m_{i+1},\dots ,m_n)\in\NN_*^n$$
such that
$$(m_1,\dots ,m_{i-1},m_{i+1},\dots ,m_n)\in\log_{\a^\prime}(b-F^ra_i),$$
where $\a^\prime=(a_1,\dots ,a_{i-1},a_{i+1},\dots ,a_n)$.
By the inductive hypothesis applied to $\a^\prime$,
$\log_{\a^\prime}(b-F^ra_i)=\bigcup_{j\in J}(\overline{r}^\prime_j\oplus B^\prime_j)$, where the $B^\prime_j\subset\NN_*^{n-1}$ come from a finite collection of basic sets that depends only on $\a^\prime$, and hence only on $\a$ and $i$.
For each $j\in J$, let $\overline{r}_j\in\NN^n$ be the tuple obtained from $\overline{r}^\prime_j$ by plugging in $0$ between the $(i-1)$st and $i$th coordinate; and let $B_j\subset\NN_*^n$ be the basic set obtained from $B^\prime_j$ by plugging in $r$ between the $(i-1)$st and $i$th coordinate of each element.
It is then not hard to see that
$$\{\overline{m}\in\log_{\a}b\colon m_i=r\}=\bigcup_{j\in J}(\overline{r}_j\oplus B_j),$$
This gives us the desired description for $\{\overline{m}\in\log_{\a}b\colon m_i=r\}$.
Since there are only finitely many choices for $i\leq t$ and $r\leq N$ (independently of $b$), this takes care of case $2$.

We are left to consider cases $3_{i,j,r}$ for fixed $i<j\leq t$ and $0\leq r\leq N$.
The set of elements of $\log_{\a}b$, for which $m_j=m_i+r$ can be described as the set of all
$$(m_1,\dots ,m_{j-1},m_i+r,m_{j+1},\dots ,m_n)\in\NN_*^n$$
such that
$$(m_1,\dots ,m_{j-1},m_{j+1},\dots ,m_n)\in\log_{\a^\prime}(b),$$
where now $\a^\prime=(a_1,\dots ,a_{i-1},a_i+F^ra_j,a_{i+1},\dots ,a_{j-1},a_{j+1},a_n)$.
The inductive hypothesis applied to $\a^\prime$ yields $\log_{\a^\prime}b=\bigcup_{l\in J}(\overline{r}^\prime_l\oplus B^\prime_l)$, where $B^\prime_l\subset\NN_*^{n-1}$ are basic sets that come from a finite collection that depends only on $\a^\prime$, and hence only on $\a ,i,j$, and $r$.
For each $l\in J$, we have two choices:
\begin{itemize}
\item[(1)] If some element of $B^\prime_l$ has nonzero $i$th coordinate then let
\begin{itemize}
\item $\overline{r}_l\in\NN^n$ be the tuple obtained from $\overline{r}^\prime_l$ by plugging in $0$ between the $(j-1)$st and $j$th coordinate;
\item $B_l\subset\NN_*^n$ be the basic set obtained from $B^\prime_l$ by plugging in the sum of the $i$th coordinate and $r$ between the $(j-1)$st and $j$th coordinate of each element.
\end{itemize}
\item[(2)] If the $i$th coordinate of every element of $B^\prime_l$ is zero then let
\begin{itemize}
\item $\overline{r}_l\in\NN^n$ be obtained from $\overline{r}^\prime_l$ by plugging in the sum of the $i$th coordinate and $r$ between the $(j-1)$st and $j$th coordinate;
\item $B_l\subset\NN_*^n$ be the basic set obtained from $B^\prime_l$ by
plugging in $0$ between the $(j-1)$st and $j$th coordinate of each element.
\end{itemize}
\end{itemize}
It is then not hard to see that
$$\{\overline{m}\in\log_{\a}b\colon m_j=m_i+r\}=\bigcup_{l\in J}(\overline{r}_l\oplus B_l),$$
and the $B_l$'s come from a finite collection of basic sets that only depend on $\a,i,j$, and $r$.
Ranging over the finitely many possibilities for $i<j\leq t$ and $0\leq r\leq N$, we complete the proof of Theorem \ref{log} (and hence of Proposition~\ref{log-closed}).
\qed

\bigskip
\bigskip

\section{Quantifier Elimination}

Fix an $F$-structure $(M,\S)$.
Our goal in this section is to prove that $\th(M,\S)$ admits quantifier elimination.
While we have already observed that $F$-sets are closed under intersections (this is part of the remark after Proposition \ref{positive-qf}), we want something more.
In order to work in elementary extensions of $M$ we will need a uniform description of such intersections.
This is where the uniformity in Theorem \ref{log} will be useful.
We begin with the special case of the intersection of a groupless $F$-set with a group:

\begin{lemma}
\label{intersect-groupless-group}
Suppose $S\subset M$ is a groupless $F$-set and $G$ is an $F$-subspace of $M$.
There is a finite sequence of groupless $F$-sets, $R_1,\dots ,R_l\subset M$, such that for any $c\in M$ there is a subset $J\subset\{1,\dots ,l\}$ and points $c_1,\dots ,c_l\in M$, such that
$$(c+S)\cap G=\bigcup_{j\in J}(c_j+R_j).$$
\end{lemma}

\proof
Taking finite unions, and translating, we may assume that $S=S(\a;\overline{\delta})$ is a basic groupless $F$-set.
We can write $S=F^{B_S}\a$, where $B_S$ is the variety $\{(m_1,\dots ,m_n)\in\NN_*^n\colon m_i\equiv 0\mod\delta_i\}$.
Let $\pi$ denote reduction modulo $G$.
Then, $(c+S)\cap G=c+F^B\a$ where $B=\log_{(-\pi\a)}\pi c\cap B_S$.
We need to describe the set $B$.

Working in $M/G$ and applying Theorem \ref{log} to $\log_{(-\pi\a)}\pi c$, we have that there exist basic sets $B_1,\dots ,B_l\subset \NN_*^{n}$ (independent of $c$), a subset $J\subset\{1,\dots ,l\}$, and tuples $\overline{r}_j$ for $j\in J$, such that
$$\log_{(-\pi\a)}\pi c=\bigcup_{j\in J}(\overline{r}_j\oplus B_j).$$
Since $B_S$ is given purely by congruence conditions, and because we are dealing with disjoint translates, it is not hard to see that for each $j\in J$,
$$(\r_j\oplus B_j)\cap B_S=\r_j\oplus(B_j\cap B_S)$$
if $\r_j\in B_S$, and is empty otherwise.
For each $i\leq l$, let $C_i=B_i\cap B_S$.
Then $C_1,\dots ,C_l$ are basic sets that are still independent of $c$.
Let $J^\prime$ be those $j\in J$ such that $\r_j\in B_S$.
We have
$$B=\log_{(-\pi\a)}\pi c\cap B_S=\bigcup_{j\in J^\prime}(\overline{r}_j\oplus C_j).$$

It follows that
$$F^B\a=\bigcup_{j\in J^\prime}(F^{(\overline{r}_j\oplus C_j)}\a)=\bigcup_{j\in J^\prime}(F^{\r_j^\circ}\a+F^{C^\circ_j}\a)$$
where $\r_j^\circ$ and $C_j^\circ$ are the canonical contractions associated to disjoint tranlsates.
Each $F^{C^\circ_j}\a$ is a basic groupless $F$-set, say $R_j$, that does not depend on $c$.
Letting $c_j=c+F^{\r_j^\circ}\a$ we obtain
$$(c+S)\cap G=\bigcup_{j\in J}(c_j+R_j),$$
as desired.
\qed

\smallskip

The following additional uniformity is inherent in the proof of the above lemma.
\begin{corollary}
\label{uniformity-groupless-group}
Let $S\subset M$ be a groupless $F$-set and $G$ an $F$-subspace of $M$.
Then there exist groupless $F$-sets $R_1,\dots ,R_l$, and for each $J\subset\{1,\dots ,l\}$ there is $U_J\subset M\times M^{|J|}$ such that
\begin{itemize}
\item[(i)] $U_J$ is a generalised $F$-set.
\item[(ii)] The {\em (}first coordinate{\em )} projections of the $U_J$'s cover $M$.
\item[(iii)] For all $(c,c_j)_{j\in J}\in U_J$,
$\displaystyle (c+S)\cap G=\bigcup_{j\in J}c_j+R_j.$
\end{itemize}
\end{corollary}

\proof
Let  $S=F^{B_S}\a$, and $\pi$ the reduction modulo $G$ map, as in the proof of Lemma~\ref{intersect-groupless-group}.
We have that for all $c\in M$, $(c+S)\cap G=c+F^B\a$, where $B$ was the closed set $\log_{(-\pi\a)}\pi c\cap B_S$.
Moreover, we had that $\displaystyle B=\bigcup_{j\in J}\r_j\oplus C_j$, the $R_j$'s were obtained as $F^{C_j^\circ}\a$, and the $c_j$'s were obtained as $c+F^{\r_j^\circ}\a$.

The main point is that the set
$$X_J=\{(\r,\s_j)_{j\in J}\colon\log_{(-\pi\a)}F^{\r}(-\pi\a) \cap B_S=\bigcup_{j\in J}\s_j\oplus C_j\}$$
is definable in $\NN_*$.
This is because $\log_{(-\pi\a)}F^{\r}(-\pi\a)$ is the $(-\pi\a)$-equivalence class of $\r$, and hence, as $\r$ varies, is uniformly definable.
Also, the $\s_j\oplus C_j$'s are uniformly definable over the $\s_j$'s.

It follows that,
$$V_J=\{(-F^{\r}\pi\a,F^{\s_j^\circ}\a)_{j\in J}\colon (\r,\s_j)_{j\in J}\in X_J\}$$
is a boolean combination of groupless $F$-sets in $M/G\times M^{|J|}$.
Consider the map from $M$ to $M/G\times M^{|J|}$ given by $(x,x_j)_{j\in J}\mapsto (x/G,x_j-x)_{j\in J}$.
Then, the preimage of $V_J$ under this map, $U_J$, is a generalised $F$-set (by Proposition \ref{positive-qf}) and satisfies the conclusion of the corollary.
\qed

\smallskip

\begin{proposition}
\label{intersect-general-group}
Suppose $S\subset M$ is a groupless $F$-set, and $H,G$ are $F$-subspaces of $M$.
Then there exist groupless $F$-sets $R_1,\dots ,R_l$, and for each $J\subset\{1,\dots ,l\}$ there is $U_J\subset M\times M^{|J|}$ such that
\begin{itemize}
\item[(i)] $U_J$ is a generalised $F$-set.
\item[(ii)] The {\em (}first coordinate{\em )} projections of the $U_J$'s cover $M$.
\item[(iii)] For all $(c,c_j)_{j\in J}\in U_J$,
$\displaystyle (c+S+H)\cap G=\bigcup_{j\in J}c_j+R_j+(H\cap G)$.
\end{itemize}
\end{proposition}

\proof
First of all, note that $(c+S+H)\cap G=[(c+S)\cap(H+G)+H]\cap G$.
Let $S_1,\dots S_l$ and $\{V_J\colon J\subset\{1,\dots ,l\}\}$ be the groupless $F$-sets and generalised $F$-sets given by Corollary \ref{uniformity-groupless-group} applied to $S$ and $H+G$.
Then for some $J$ and all $(c,d_j)_{j\in J}\in V_J$ we have
$$(c+S)\cap(H+G)=\bigcup_{j\in J}d_j+S_j$$
where now, $d_j$ and $S_j$ are in $H+G$.
Let $\pi$ denote reduction modulo $H\cap G$.
Fixing $j\in J$, note that
$$(d_j+S_j+H)\cap G=\pi^{-1}[(\pi d_j+\pi S_j+\pi H)\cap \pi G].$$
But as $\pi H+\pi G=\pi H\oplus\pi G$, letting $\eta\colon\pi H+\pi G\to\pi G$ be the direct sum coordinate projection, we see that
$$(\pi d_j+\pi S_j+\pi H)\cap \pi G=\eta\pi d_j+\eta\pi S_j.$$
Now, $\pi^{-1}(\eta\pi S_j)=R_j+H\cap G$, where $R_j$ is a groupless $F$-set.
It follows that
$$(d_j+S_j+H)\cap G=c_j+R_j+H\cap G,$$
and this holds for any $(d_j^\prime,c_j^\prime)\in (H+G)^2$ such that $\pi c_j^\prime=\eta\pi d_j^\prime$.
Putting all this together we obtain
\begin{eqnarray*}
(c+S+H)\cap G
& = & [(c+S)\cap(H+G)+H]\cap G \\
& = & \bigcup_{j\in J}[(d_j+S_j+H)\cap G] \\
& = & \bigcup_{j\in J}c_j+R_j+(H\cap G)
\end{eqnarray*}
Let $U_J$ be the set of all $(x,x_j)_{j\in J}$ such that for some $(y_j)_{j\in J}$, $(x,y_j)_{j\in J}\in V_J$ and $\pi x_j=\eta\pi y_j$.
Then $U_J$ is also a generalised $F$-set, and the above equality holds for all $(c,c_j)_{j\in J}\in U_J$.
This completes the proof of the proposition.
\qed

\smallskip

\begin{corollary}
\label{intersect-general}
Suppose $S,T\subset M$ are groupless $F$-sets, and $H,G\subset M$ are $F$-subspaces.
Then there exist groupless $F$-sets $R_1,\dots ,R_l$ such that for all $c\in M$ there is $J\subset\{1,\dots ,l\}$ and points $(c_j)_{j\in J}$ from $M$, such that
$$(c+S+H)\cap(T+G)=\bigcup_{j\in J}c_j+R_j+(H\cap G).$$
\end{corollary}

\proof
This follows from the main conclusion of Proposition \ref{intersect-general-group} (ignoring the $U_J$'s) exactly as in case $3$ of Proposition \ref{positive-qf}.
That is, as $(c+S+H)\times(T+G)$ is again an $F$-set, and as the diaogonal $\Delta\subset M^2$ is an $F$-suubspace, we get a uniform description of $[(c+S+H)\times(T+G)]\cap\Delta$ by Proposition \ref{intersect-general-group}.
Projecting back to the first coordinate, we get a uniform description for $(c+S+H)\cap(T+G)$, as desired.
\qed

\smallskip

\begin{remark}
\label{intersect-nonstandard}
If ${^*}M$ is an elementary extension of $M$, then by an $F$-set (or $F$-subspace) of ${^*}M$ we mean the interpretation in ${^*}M$ of an $F$-set (or $F$-subspace) in $M$.
The above corollary implies that if $X,Y\subset{^*}M$ are finite unions of translates of $F$-sets, then so is $X\cap Y$.
Moreover, the $F$-subspaces appearing in $X\cap Y$ are the intersections of $F$-subspaces appearing in $X$ and in $Y$.
\end{remark}

\bigskip

We have one other ingredient before we can address quantifier elimination.
We need to be able to describe when a given $F$-set is covered by a union of translates of given $F$-sets.
We begin, as usual, with the groupless case.

\begin{proposition}
\label{cover-translates}
Suppose $T,S_1,\dots ,S_l\subset M$ are groupless $F$-sets.
Then the set
$$U=\{(c_1,\dots c_l)\in M^l\colon T\subset\bigcup_{i=1}^l(c_i+S_i)\}$$
is a boolean combination of groupless $F$-sets.
\end{proposition}

\proof
There are closed sets $B\subset\NN_*^p$ and $C_1,\dots ,C_l\subset\NN_*^r$, and tuples $\x\in M^p$ and $\z\in M^r$, such that $T=F^B\x$ and $S_i=F^{C_i}\z$.
Let $U^\prime$ be those $(c_1,\dots ,c_l)\in U$, such that $c_i+S_i$ has nonempty intersection with $T$ for all $i$ -- that is, there are no extraneous sets involved.

Note that if $(c_1,\dots ,c_l)\in U^\prime$, then $c_i=F^{\r_i}\x-F^{\s_i}\z$ for some tuples $\r_i\in B$ and $\s_i\in C_i$.
Letting $\y=(\x,-\z)$ and $q=p+r$, we have that each $c_i$ is of the form $F^{\r_i}\y$ for some $r_i\in\NN_*^q$.
That is, letting
$$V=\{(F^{\r_1}\y,\dots ,F^{\r_l}\y)\colon F^B\x\subset\bigcup_{i=1}^lF^{(\r_i\times C_i)}(\y,\z)\},$$
we obtain $U^\prime\subset V\subset U$.
It is sufficient to prove that $V$ is a finite boolean combination of groupless $F$-sets, since $U$ will be a finite union of such $U^\prime$'s -- and hence of such $V$'s -- as we vary among the subsets of $\{1,\dots ,l\}$.

Let $\w=(\x,\y,\z)$, and let $A\subset (\NN_*^q)^l$ be those tuples $(\r_1,\dots ,\r_l)$ such that
$$B\times\{-\infty\}^{q+r}\subset\bigcup_{i=1}^l(\{-\infty\}^p\times\r_i\times C_i) \mod\sim_{\w}.$$
Then,
$$V=\{(F^{\r_1}\y,\dots ,F^{\r_l}\y)\colon (\r_1,\dots ,\r_l)\in A\}.$$
Since $A$ is definable, $V$ is a finite boolean combination of groupless $F$-sets.
\qed

\begin{corollary}
\label{cover-translates-cor}
Suppose $T,S_1,\dots ,S_l\subset M$ are groupless $F$-sets.
Then the set
$$V=\{(c_0,c_1,\dots c_l)\in M^{l+1}\colon c_0+T\subset\bigcup_{i=1}^l(c_i+S_i)\}$$
is a generalised $F$-set.
\end{corollary}

\proof
Clearly,
$$c_0+T\subset\bigcup_{i=1}^l(c_i+S_i)\iff T\subset\bigcup_{i=1}^l[(c_i-c_0)+S_i]$$
Hence, if we let $U$ be the boolean combination of groupless $F$-sets given by Proposition \ref{cover-translates} applied to $T,S_1,\dots S_l$, and $\eta\colon M^{l+1}\to M^l$ the surjection given by $\eta(a_0,\dots a_l)=(a_1-a_0,\dots ,a_l-a_0)$; then $V=\eta^{-1}(U)$.
This is a generalised $F$-set by Proposition \ref{positive-qf}.
\qed

\medskip

We want to prove Corollary \ref{cover-translates-cor} with arbitrary $F$-sets replacing groupless $F$-sets.
We need several preliminary lemmas.

\begin{lemma}
\label{nofgroup}
If $G$ is an infinite subgroup of $M$, then $G$ is \emph{not} a groupless $F$-set.
\end{lemma}

\proof
Suppose $G$ is a groupless $F$-sets.
Then for some $\a\in\NN_*^n$ and closed $B\subset\NN_*^n$, $G=F^B\a$.
Consider the relation $P\subset B\times B\times B$ given by $(\r_1,\r_2,\r_3)\in P$ if and only if $F^{\r_1}\a+F^{\r_2}\a=F^{\r_3}\a$.
It is not hard to see that
$$P=\{(\r_1,\r_2,\r_3)\colon (\r_1,\r_2)\sim_{(\a,\a)}(\r_3,-\infty)\}.$$
Since $(\a,\a)$-equivalence is definable in $\NN_*$, $P$ is a definable set.
Now $P$ defines a group structure on $B/\sim_{(\a,\a)}$ that is isomorphic to that of $G$.
However, $\NN_*$ does not interpret an infinite group.
\qed

\smallskip

\begin{lemma}
\label{finind}
Suppose $H_1, \dots, H_n\leq G$ are $F$-subspaces of $M$ such that each $H_i$ is of infinite index in $G$, 
and $S_1, \ldots, S_n\subset M$ are groupless $F$-sets.
Then
$$G \neq\bigcup_{i=1}^n S_i + H_i.$$
\end{lemma}

\proof
We proceed by induction on $n$.
Consider the case of $n=1$.
Let $\pi$ be the reduction modulo $H_1$ map.
If $G=S_1+H_1$, then $\pi G=\pi S_1$.
But as $\pi G$ is infinite, this would contradict Lemma \ref{nofgroup}.

Suppose $n>1$.
Define a quasiodering on ${1, \ldots, n}$ by 
$i \sqsubseteq j$ if and only if $H_i \cap H_j$ has finite index in $H_i$.
{\em We may assume that $\sqsubseteq$ is a partial ordering.}
If not, then it must be the case that for some $i\neq j$, $H_i\cap H_j$ is of finite index in both $H_i$ and $H_j$.
But then, $(S_i+H_i)\cup (S_j+H_j)=S^\prime+(H_i\cap H_j)$ for some groupless $F$-set $S^\prime$.
But by induction, $\displaystyle G \neq [S^\prime+(H_i\cap H_j)]\cup \bigcup_{k\neq i,j} S_k + H_k$.

After reordering, we may assume that $H_1$ is $\sqsubseteq$-maximal.
That is, for all $i\neq 1$, $H_i\cap H_1$ is of infinite index in $H_1$.
Let $\pi\colon G\to G/H_1$ be the quotient map restricted to $G$.
Now, as $H_1$ has infinite index in $G$, $X=\pi(S_1\cup H_1)=\pi S_1\neq G/H_1$ by Lemma \ref{nofgroup}.
Let $h\in (G/H_1)\setminus X$.
Then $\pi^{-1}(h)$ is disjoint from $S_1+H_1$.
Hence, if $\displaystyle G = \bigcup_{i=1}^n S_i + H_i$ then $\pi^{-1}(h)$ is covered by $\displaystyle \bigcup_{i=2}^n S_i + H_i$.
But $\pi^{-1}\{h\}$ is a coset of $H_1$, say $c+H_1$.
So $\displaystyle H_1=\bigcup_{i=2}^n (c+ S_i + H_i)\cap H_1$.
Expressing the $(c+ S_i + H_i)\cap H_1$'s as $F$-sets using Lemma \ref{intersect-general-group}, and recalling that for $i\neq 1$ $H_1\cap H_i$ has infinite index in $H_1$; we see that this is impossible by induction.
\qed

\smallskip

\begin{corollary}
\label{repfin}
Suppose $H_1, \dots, H_n\leq G$ are $F$-subspaces and $S_0, \ldots, S_n \subset G$ are groupless $F$-sets.
Let $I$ be the set of indices $i\leq n$ such that $H_i$ has finite index in $G$.
If $\displaystyle G=S_0\cup \bigcup_{i=1}^n(S_i + H_i)$,
then $\displaystyle G=\bigcup_{i \in I}(S_i + H_i)$.
\end{corollary}

\proof
Let $\displaystyle K=\bigcap_{i\in I}H_i$, and $J=\displaystyle \bigcup_{i\in I} S_i + H_i$.
Then $J$ is a finite union of cosets of $K$.
Suppose that $G\setminus J\neq\emptyset$, and $a \in G \setminus J$.
We would have $a + K \subseteq S_0 \cup \bigcup_{i \notin I} S_i + H_i$, so that $K = [(S_0 - a) \cap K] \cup \bigcup_{i \notin I} [(S_i - a) + H_i] \cap K$ contradicting Lemma~\ref{finind}.
\qed

\smallskip

We are now able to extend Corollary \ref{cover-translates-cor} to arbitrary $F$-sets.

\begin{proposition}
\label{cover-translates-general}
Suppose $G,H_1,\dots ,H_l$ are $F$-subspaces of $M$ and $T, S_1,\dots ,S_l$ are groupless $F$-sets.
Then the set
$$V=\{(c_0,c_1,\dots c_l)\in M^{l+1}\colon c_0+T+G\subset\bigcup_{i=1}^l(c_i+S_i+H_i)\}$$
is a generalised $F$-set.
\end{proposition}

\proof
Clearly $c_0+T+G$ is covered by the union on the righthand side of the displayed inclusion if and only if $c_0+t+G$ is covered for all $t\in T$.
But by Corollary~\ref{repfin}, if this happens at all, it will already happen with those indices $i$ such that $H_i \cap G$ has finite index in $G$.
Thus, we may and do assume that $H_i \cap G$ is of finite index in $G$ for each $i$.
Let $\displaystyle K := \bigcap_{i=1}^l (H_i \cap G)$.
Denote the reduction modulo $K$ by $\pi$.
Then $c_0+T+G\subset\bigcup_{i=1}^l(c_i+S_i+H_i)$ if and only if $\pi c_0+\pi T+\pi G\subset\bigcup_{i=1}^l(\pi c_i+\pi S_i+\pi H_i)$.
Since $\pi G$ and the $\pi H_i$'s are finite sets, we are in the groupless case.
By Corollary \ref{cover-translates-cor}, the set
$$U=\{(d_0,d_1,\dots ,d_l)\colon d_0+\pi T+\pi G\subset\bigcup_{i=1}^l(d_i+\pi S_i+\pi H_i)\}\subset (M/K)^{l+1}$$
is a generalised $F$-set.
Hence,
$$V=\{(c_0,c_1,\dots c_l)\in M^{l+1}\colon c_0+T+G\subset\bigcup_{i=1}^l(c_i+S_i+H_i)\}=\pi^{-1}(U)$$
is also a generalised $F$-set.
\qed

\smallskip

\begin{remark}
\label{abelian-structure}
Note that if the $F$-subspaces that appear in the hypotheses of Propositions \ref{intersect-general-group} and \ref{cover-translates-general} come from a class of $F$-subspaces $\A$ such that $(M,\A)$ is an abelian structure admitting quantifier elimination, then the $F$-subspaces that appear in the generalised $F$-sets of their respective conclusion also come from $\A$.\footnote{See the discussion at the end of Section \ref{set-up}.}
\end{remark}

\bigskip

We have described, uniformly and in terms of generalised $F$-sets,
\begin{itemize}
\item how to write the intersection of a translate of an $F$-set with an $F$-subspace as a union of translates of $F$-sets (this is Proposition \ref{intersect-general-group}), and
\item when a translate of an $F$-set is covered by a union of translates of given $F$-sets (this is Proposition \ref{cover-translates-general}).
\end{itemize}
We are now in a position to put these together to obtain:

\begin{proposition}
\label{cover-intersect}
Given groupless $F$-sets $T,S_1,\dots S_l\subset M$, and $F$-subspaces $L,G,H_1,\dots ,H_l$,
$$V=\{c\in M\colon(c+T+G)\cap L\subset\bigcup_{i=1}^l[(c+S_i+H_i)\cap L]\}$$
is a generalised $F$-set.
\end{proposition}

\proof
Let $T_1,\dots T_n$ be the groupless $F$-sets given by Proposition \ref{intersect-general-group} applied to $T$, $G$, and $L$;
and let $R_{i1},\dots R_{im_i}$ be the groupless $F$-sets given by that Proposition applied to $S_i$, $H_i$, and $L$, for each $i\leq l$.
Then $c\in V$ if and only if for some $J\subset\{1,\dots n\}$, and $J_i\subset\{1,\dots m_i\}$,
$$(c+T+G)\cap L=\bigcup_{j\in J}(a_j+T_j+G\cap L),$$
and 
$$(c+S_i+H_i)\cap L=\bigcup_{j\in J_i}(a_{ij}+R_{ij}+H_i\cap L),$$
for each $i\leq l$, and
$$\bigcup_{j\in J}(a_j+T_j+G\cap L)\subset\bigcup_{i=1}^l\bigcup_{j\in J_i}(a_{ij}+R_{ij}+H_i\cap L).$$
By Proposition \ref{intersect-general-group} there are generalised $F$-sets $U_J$ and $U_{J_i}$ such that the first two equalities hold for all $(c,a_j)_{j\in J}\in U_J$ and $(c,a_{ij})_{j\in J_i}\in U_{J_i}$.
By Proposition \ref{cover-translates-general} for fixed $j^\prime\in J$, the set of $(a_{j^\prime},a_{ij})_{i\leq l, j\in J_i}$'s such that the final displayed containment holds is also a generalised $F$-set.
As generalised $F$-sets are closed under intersections, unions, and projections, $V$ is a generalised $F$-set.
\qed

\bigskip

\begin{theorem}
\label{qe}
The theory of an $F$-structure admits quantifier elimination.
\end{theorem}

\proof
Suppose $(M,\S)$ is an $F$-structure.
We need to show that the projection of a boolean combination of $F$-sets is again a boolean combination of $F$-sets.
As $F$-sets are closed under intersections and union, and as projections commute with intersections and unions, it is sufficient to show that
$\displaystyle \pi[(T+G)\setminus\bigcup_{i=1}^l(S_i+H_i)]$
is a boolean combination of $F$-sets; where $T, S_1,\dots S_l\subset M^{m+n}$ are groupless $F$-sets, $G,H_1,\dots H_l\leq M^{m+n}$ are $F$-subspaces, and $\pi\colon M^{m+n}\to M^m$ is the projection map.
Now,
$$\pi[(T+G)\setminus\bigcup_{i=1}^l(S_i+H_i)]=\pi(T+G)\setminus V$$
where
$$V=\{c\in M^m\colon (T+G)_c\subset\bigcup_{i=1}^l(S_i+H_i)_c\}.$$
Hence, the Theorem will follow once we know that $V$ is a boolean combination of $F$-sets.
Notice that $c\in V$ if and only if
$$(T+G)\cap(c\times M^n)\subset\bigcup_{i=1}^l(S_i+H_i)\cap(c\times M^n)$$
if and only if
$$[(-c,0)+T+G]\cap(0\times M^n)\subset\bigcup_{i=1}^l[(-c,0)+S_i+H_i]\cap(0\times M^n).$$
It follows from Proposition \ref{cover-intersect} that the set of $(-c,0)$'s for which this occurs is a generalised $F$-set.
Hence, $V$ is a generalised $F$-set, and so, in particular, is a boolean combination of $F$-sets (see Remark \ref{generalised-qf}).
This completes the proof of quantifier elimination.
\qed

\smallskip

\begin{remark}
\label{abelian-structure-qe}
If $\A$ is a collection of $F$-subspaces of cartesian powers of $M$ such that $(M,\A)$ is an abelian structure admitting quantifier elimination, then $\th(M,\S\A)$ -- that is, the theory of the $F$-structure on $(M,\A)$ -- also admits quantifier elimination.
Indeed, the quantifier elimination process presented here does not produce any new $F$-subspaces (see Remark \ref{abelian-structure}).
\end{remark}

\bigskip
\bigskip

\section{Stability}

The goal of this section is to use quantifier elimination to show that the theory of an $F$-structure is stable.
Fix an $F$-space $M$.
We introduce a stratification by complexity of the groupless $F$-sets of $M$:

\begin{definition}
Suppose $\delta\in\NN$ with $\delta>0$.
We define the class of groupless $F$-sets in $M$ {\em of order $\delta$}, denoted by $O(\delta)$ to be those $F$-sets that are obtained as finite unions of translates of basic groupless $F$-sets of the form $S(a_1,\dots ,a_n;\delta_1,\dots ,\delta_n)$, where each $\delta_i$ divides $\delta$.
\end{definition}

The point here is that any groupless $F$-set in $O(\delta)$ can be written as a finite union of groupless $F$-sets that only involve $\delta$.
The relationship between $\L^\delta$-definable sets and groupless $F$-sets is enhanced if we restrict our attention to fixed strata.
Let us fix some further notation that is convenient for discussing this relationship.
For any set $Y$ from $\NN_*$ and any tuple $\a$ from $M$, we will use $[Y]_{\a}$ to denote the quotient of $Y$ by $\a$-equivalence.

\begin{definition}
If $S \subset M$ is any groupless $F$-set and $\a$ is an $n$-tuple from $M$, then we say that {\em $\a$ sees $S$} if $S=F^B\a$ for some closed $B\subset\NN_*^n$.
In this case we use the notation:
$$\log_{\a}S=\{\r\in\NN_*^n\colon F^{\m}\a\in S\}.$$
Note that $\log_{\a}S$ is just the saturation of $B$ by $\a$-equivalence, and (as $\a$-equivalence is a closed relation) is also a witness to the fact that $\a$ sees $S$.
That is, $\log_{\a}S$ is closed and $[\log_{\a}S]_{\a}=[B]_{\a}$.
\end{definition}

\smallskip

Recall that for any tuple $\a$ from $M$, $\a$-equivalence is a definable relation in $\NN_*$ (this is Corollary \ref{equivalence} and follows from Theorem \ref{log}).
In fact, if $\delta_0>0$ is such that $F^{\delta_0}$ fixes $F^\infty M$, then for all tuples $\a$ from $M$, $\a$-equivalence is a $\delta_0$-closed relation.\footnote{Congruence conditions only entered into the description of $\log_{\a}b$ in Theorem~\ref{log} in the base case dealt with by Lemma \ref{special-log} -- where one worked with a power of $F$ that fixes $F^\infty M$.}
That is, these equivalence relations are all definable in the fragment $\L^{\delta_0}$ of $\L^*$.
For the rest of this section we fix such a $\delta_0$.

\begin{lemma}
\label{well-defined-dimension}
Suppose $\delta$ is a multiple of $\delta_0$, and $S\in O(\delta)$.
There is a tuple $\a$ from $M$ that sees $S$, such that $\log_{\a}S$ is $\L^\delta$-definable.
Moreover, if $\b$ is any other tuple that sees $S$, then $\log_{\b}S$ is $\L^\delta$-definable and there is an $\L^\delta$-definable bijection between $[\log_{\a}S]_{\a}$ and $[\log_{\b}S]_{\b}$.
\end{lemma}

\proof
By Lemma \ref{closed-exponent=groupless}, every groupless $F$-set is seen by some tuple from $M$.
Hence $S=F^B\a$ for some $\a\in M^m$ and some closed $B\subset\NN_*^m$.
It is not hard to see (from the proof of Lemma \ref{closed-exponent=groupless}, for example), that $\a$ and $B$ can be chosen so that $B$ is in fact $\delta$-closed.
Since $\delta_0$ divides $\delta$, $\a$-equivalence is $\L^\delta$-definable.
Hence, $\log_{\a}S$, being the saturation of $B$ with respect to $\a$-equivalence, is also $\L^\delta$-definable.

Now suppose that $\b\in M^n$ also sees $S$.
We need to show that $\log_{\b}S$ is $\L^\delta$-definable and that there is an $\L^\delta$-definable bijection between $[\log_{\a}S]_{\a}$ and $[\log_{\b}S]_{\b}$.
Consider the relation
$$X=\{(\r,\s)\colon\r\in\NN_*^m, \s\in\NN_*^n, F^{\r}\a=F^{\s}\b\}.$$
We claim that $X$ is $\L^\delta$-definable.
Indeed, let $\c\in M^{m+n}$ be the concatenation of $\a$ and $\b$.
Now $(\r,\s)\in X$ if and only if $F^{(\r,-\infty)}\c=F^{(-\infty,\s)}\c$.
But $\c$-equivalence is $\L^\delta$-definable.

Observe that $\s\in\log_{\b}S$ if and only if for some $\r\in\log_{\a}S$, $(\r,\s)\in X$.
Since $\log_{\a}S$ and $X$ are $\L^\delta$-definable, this shows that $\log_{\b}S$ is also $\L^\delta$-definable.
Finally, it is not hard to see that $X$ descends and restricts to the graph of a bijection between $[\log_{\a}S]_{\a}$ and $[\log_{\b}S]_{\b}$.
Hence $[\log_{\a}S]_{\a}$ and $[\log_{\b}S]_{\b}$ are $\L^\delta$-definably isomorphic.
\qed

\smallskip

For groupless $F$-sets restricted to a given stratum, we can now use Morley rank in $\NN_*$ to define a dimension.
The following notions of $\delta$-dimension and $\delta$-degree are well defined by Lemma \ref{well-defined-dimension}

\begin{definition}
Fix $\delta$ a multiple of $\delta_0$.
Let $S \subset M$ be a groupless $F$-set in $O(\delta)$.
The {\em $\delta$-dimension of $S$} (respectively,  {\em $\delta$-degree of $S$}), denoted by $\dim_{\delta}S$  (respectively $\deg_{\delta}S$), is the Morley rank (respectively Morley degree) of $[\log_{\a}S]_{\a}$, in the sense of $\L^\delta$, for any $\a$ that sees $S$.
\end{definition}

\smallskip

That this notion of dimension is well behaved, even in elementary extensions, is exemplified in the following proposition:

\begin{proposition}
\label{closure-congruence-free}
Let $\delta$ be a multiple of $\delta_0$.
Suppose $S$ and $T$ are groupless $F$-sets in $O(\delta)$ with $\dim_\delta S=\dim_\delta T$.
Then there are groupless $F$-sets $R_1,\dots ,R_l$ in $O(\delta)$, such that for all $i$, $\dim_\delta R_i<\dim_\delta T$ or $\deg_\delta R_i<\deg_\delta T$; and such that the following holds:  for all $c\in M$, if $c+S\subset T$ then there is $I\subset\{1,\dots ,l\}$ and $(b_i)_{i\in I}$ such that $\displaystyle [T\setminus (c+S)]\subset\bigcup_{i\in I} b_i+R_i$.
\end{proposition}

\proof
We first claim that there is a unifomly definable family of $\L^{\delta}$-definable sets, $\{Y_{\r}\subset\NN_*^n\colon\r\in\NN_*^m\}$, such that if $c\in M$ with $c+S\subset T$, then for some $\r$, $T\setminus(c+S)=F^{Y_{\r}}(\a)$, and either the Morley rank or Morley degree of $Y_{\r}$ is strictly less than $\dim_\delta T$ or $\deg_\delta T$ (respectively).
Indeed, write $T=F^B\a$, where $\a\in M^n$ and $B\subset\NN_*^n$ is $\delta$-closed.
We may assume that $B=\log_{\a}T$, that is, $B$ is saturated for $\a$-equivalence.
Using arguments by now familiar, there is a uniformly definable family of $\delta$-closed sets, $\{D_{\r}\subset\NN_*^n\colon\r\in\NN_*^m\}$, such that if $c\in M$ with $c+S\subset T$, then for some $\r$, $c+S=F^{D_{\r}}(\a)$.
We may assume that the $D_{\r}$ are also saturated for $\a$-equivalence.
Hence, $[B]_{\a}\setminus [D_{\r}]_{\a}=[B\setminus D_{\r}]_{\a}$.
Letting $X_{\r}=B\setminus D_{\r}$, it folows that either the Morley rank or Morley degree of $[X_{\r}]_{\a}$ is strictly less than $\dim_\delta T$ or $\deg_\delta T$, respectively.
Moreover, $T\setminus(c+S)=F^{X_{\r}}(\a)$.
Using the existence of definable Skolem functions in $(\NN_*,\sigma,0,-\infty,P_\delta)$, we find definable $Y_{\r}\subset X_{\r}$ such that $\rm(Y_{\r})=\rm[X_{\r}]_{\a}$, $\dm(Y_{\r})=\dm[X_{\r}]_{\a}$, and $T\setminus(c+S)=F^{Y_{\r}}(\a)$.

Using part $(d)$ of Lemma \ref{exponent-qe}, we obtain a finite union of $\delta$-varieties $Z\subset\NN_*^{m+n}$ with $Y_{\r}\subset Z_{\r}$, $\rm(Z_{\r})=\rm (Y_{\r})$, and $\dm(Z_{\r})=\dm(Y_{\r})$.
It follows that for each $c$ such that $c+S\subset T$, there is an $\r$ such that $T\setminus(c+S)\subset F^{Z_{\r}}(\a)$ and either $\rm Z_{\r}<\dim_\delta T$ or $\dm Z_{\r}<\deg_\delta T$.
By Lemma \ref{uniform-variety=disjoint}, $\{Z_{\r}\}$ is a uniform family of disjoint translates of a fixed $\delta$-closed set.
It follows that $\{F^{Z_{\r}}(\a)\}$ is a uniform family of unions of translates of fixed groupless $F$-sets in $O(\delta)$.
These groupless $F$-sets will have either $\delta$-dimension or $\delta$-degree strictly less than that of $T$.
\qed

\bigskip

We are now in a position to prove the stability of $F$-structures.

\begin{theorem}
\label{stable}
Let $(M,\S)$ be an $F$-structure
Let $N \succeq M$ be an elementary extension as $F$-structures.
If $p(x) \in S(N)$ is a type over $N$, then $p$ is $N$-definable.
That is, the theory of $(M,\S)$ is stable.
\end{theorem}

\proof
By quantifier elimination, every definable set in $N$ is a boolean combination of fibres of $F$-sets from the ground model.
Such fibres can be expressed as intersections of translates of $F$-sets.
It follows that every definable set in $N$ is a boolean combination of translates of $F$-sets (see Remark \ref{intersect-nonstandard}).

Let $S$ be a groupless $F$-set and $G$ an $F$-subspace (in the sense of the ground model $(M,\S)$).
Consider the set 
$$D_{S,G}(p) := \{ b \in N : p(x) \vdash  x \in b + S + G \}.$$
By quantifier elimination,
it suffices to show that $D_{S,G}(p)$ is $N$-definable for all such $S$ and $G$.  
If $p = \tp(\alpha/N)$, then we set $p/G := \tp([\alpha/G]/N)$.
Let $\pi$ denote reduction modulo $G$.
If we know that $D_{\pi(S),0}(p/G)$ is $N/G$-definable, then $D_{S,G}(p)$ is $N$-definable.
Thus it suffices to show that for any groupless $F$-set $S$ and any type $p \in S(N)$ the set
$$D_S(p):= \{ b \in N : p(x) \vdash  x \in b + S\}.$$
is $N$-definable.
Clearly, it is sufficient to show this for every groupless $F$-set in $O(\delta)$, for all $\delta$ a multiple of $\delta_0$.
Fix such a $\delta$.

If for all $T$ in $O(\delta)$ the sets $D_T(p)$ are empty,
then they are certainly $N$-definable.
So, we may assume that some $D_T(p) \neq \varnothing$.
Let $T\in O(\delta)$ have minimal $\delta$-dimension and then minimal $\delta$-degree, such that for some $b\in N$, $p(x) \vdash  x \in b + T$. 
Fix such a choice of $b$.
We will prove that for all $S\in O(\delta)$, $D_S(p)$ is $b$-definable, and hence $N$-definable as desired.

Let $S$ be any groupless $F$-set in $O(\delta)$.
Let $R_1, \ldots, R_m$ be a groupless $F$-sets in $O(\delta)$ such that for any $x$ and $y$ we have $\displaystyle (x + S) \cap (y+T) = \bigcup_{i \in I} x_i + R_i$ for some $I \subseteq \{1, \ldots, m \}$ and $x_1, \ldots, x_m$.
Let $J := \{ i \leq m : \dim_\delta R_i = \dim_\delta T \}$.
We claim that
$$D_S(p) = \{ a : (\exists y) \bigvee_{j \in J} y + R_j \subseteq (a + S) \cap (b + T) \}.$$ 
This would prove that $D_S(p)$ is $b$-definable, as desired. 
For one direction, suppose that $a \in D_S(p)$.
Than $p(x) \vdash x \in (a + S) \cap (b + T)$.  
We write
$$(a + S) \cap (b + T) = \bigcup_{i \in I} c_i + R_i$$
for some $c_1, \ldots, c_m \in N$ and $I \subseteq \{1, \ldots, m \}$.
As $p(x)$ is consistent, $I \neq \varnothing$.
As $p$ is complete, $p(x) \vdash x \in c_i + R_i$ for some $i$.
By minimality of $\dim_\delta T$, $\dim_\delta R_i = \dim_\delta T$.

Conversely, suppose that $c + R_j \subseteq (a + S) \cap (b + T)$, for some some $c \in N$ and some $j\in J$.
If $a \notin D_S(p)$, then $p(x) \vdash x \in b+ (T \setminus [(c-b) + R_j])$.  By Proposition~\ref{closure-congruence-free}, there
are $d_1, \ldots, d_\ell \in N$ and groupless $F$-sets $U_1, \ldots, U_\ell$ in $O(\delta)$, such that 
$\displaystyle (T \setminus [(c-b) + R_j])\subset\bigcup_{i=1}^{\ell}d_i+U_i$, and for each $i\leq\ell$, either $\dim_\delta U_i<\dim_\delta T$ or $\deg_\delta U_i<\deg_\delta T$.
Indeed, Proposition~\ref{closure-congruence-free} tells us this for the ground model in a uniform manner, and hence it holds in the elementary extension.
Now for some $i\leq\ell$, $p(x) \vdash x \in d_i+U_i$, contradiciting the minimal choice of $T$.
This completes the proof of the theorem.
\qed

\bigskip
\bigskip
\section{Mordell-Lang for Isotrivial Semiabelian Varieties}
\label{intended-application}

In this section we prove a version of the absolute Mordell-Lang conjecture for isotrivial semiabelian varieties in terms of $F$-sets.
Using the theory of $F$-structures developed in the previous sections, we can deduce from this a description of the definable sets in the structure induced on a finitely generated $F$-submodule of a semiabelian variety over a finite field.
In particular we obtain stability, which in turn implies uniformity in the Mordell-Lang statement.

\smallskip

Recall the context of Example \ref{intended}:
$G$ is a semiabelian variety over the finite field ${\mathbb F}_q$ of characteristic $p>0$, $F:G \to G$ is the algebraic endomorphism induced from the $q$-power Frobenius, and $R=\ZZ[F]$ is the subring of the endomorphism ring of $G$ generated by $F$.
Let $K$ be a finitely generated regular extension field of ${\mathbb F}_q$, $\Gamma \leq G(K)$ a finitely generated $R$-submodule.
Recall that $\Gamma$ is an $F$ space and and every $R$-submodule of $\Gamma$ is an $F$-subspace.
We consider the $F$-structure $(\Gamma,\S)$.
Our first goal is to show that the intersection of $\Gamma^n$ with a subvariety of $G^n$ is an $F$-set in $\S$.
Indeed, we will show that these intersections are $F$-sets of an even simpler form.
We begin with several preliminary lemmas about the nature of $F$-sets in this particular context.

\begin{lemma}
\label{zc-basic-f-set}
The Zariski closure of a basic $F$-set in $\Gamma$ is defined over $\FF_q$.
\end{lemma}

\proof
Suppose $U=S(\a;\overline\delta)+H_{\Gamma}$, where $H_{\Gamma}$ is an $F$-submodule of $\Gamma$, $\a=(a_1,\dots a_n)$ is a tuple of points from $\Gamma$,
and  $\overline\delta=(\delta_1,\dots ,\delta_n)\in\NN^n$.
Let $\delta=\delta_1\cdots\delta_n$.
For each $i>0$, let $U_i=F^{i\delta}U$, and let $X_i$ be the Zariski closure of $U_i$.
Then $U_{i+1}\subset U_i$, and so $X_{i+1}\subset X_i$.
Hence, for some $N>0$, $X_{N+1}=X_N$.
But then,
$$F^{\delta}X_N=\overline{F^{\delta}U_N}=\overline{U_{N+1}}=X_{N+1}=X_N.$$
That is, $X_N$ is fixed by $F^{\delta}$, and so is defined over a finite field.
Moreover,
$$\overline{U}=\overline{F^{-N\delta}U_N}=F^{-N\delta}X_N=X_N.$$
Hence $\overline{U}$ is defined over a finite field.
But, as it is the Zariski closure of some points of $G(K)$, it is also defined over $K$.
It follows that $\overline{U}$ is defined over $K\cap\FF_q^{\alg}=\FF_q$, as desired.
\qed

\begin{lemma}
\label{congruence-free-zc}
Suppose $U=b+S(\a;\overline\delta)+H_{\Gamma}$ is a translate of a basic $F$-set in $\Gamma$.
Let $H\leq G$ be the Zariski closure of $H_{\Gamma}$.
Then $b+S(\a;1)+(H\cap\Gamma)$ is contained in the Zariski closure of $U$.
\end{lemma}

\proof
We proceed by induction on the arity of $\a$.
Suppose $U=b+S(a;\delta)+H_{\Gamma}$.
For any $c\in S(a;1)$ there is some $r\leq 0$ such that $F^{r}c\in S(a;\delta)$.
As $H_{\Gamma}$ is closed under $F$, $F^r(c+H_{\Gamma})\subset S(a;\delta)+H_{\Gamma}$.
Let $Y$ be the Zariski closure of $S(a;\delta)+H_{\Gamma}$.
It follows that $F^r\overline{(c+H_{\Gamma})}\subset Y$.
As $Y$ is defined over $\FF_q$ (Lemma \ref{zc-basic-f-set}), 
$c+H=\overline{(c+H_{\Gamma})}\subset Y$, and so $c+(H\cap\Gamma)\subset Y$.
We have shown that $c+(H\cap\Gamma)\subset \overline{S(a;\delta)+H_{\Gamma}}$, for all $c$ in $S(a;1)$.
This in turn implies that
$$b+S(a;1)+(H\cap\Gamma)\subset b+\overline{S(a;\delta)+H_{\Gamma}}=\overline{U}.$$

Now suppose that $U=b+S(a_0,\a;\delta_0,\overline\delta)+H_{\Gamma}$, and let $X$ be the Zariski closure of $U$.
Note that $U=b+S(a_0;\delta_0)+S(\a;\overline\delta)+H_{\Gamma}$.
For all $c\in S(a_0;\delta_0)$,
$$b+c+S(\a;\overline\delta)+H_{\Gamma}\subset X.$$
By induction, $b+c+S(\a;1)+(H\cap\Gamma)\subset X$.
Hence,
$$b+S(a_0;\delta_0)+S(\a;1)+(H\cap\Gamma)\subset X.$$
Similarly, for all $d\in S(\a;1)$, $b+S(a_0;\delta_0)+d+H_{\Gamma}$ is contained in $X$ -- and so by induction again, $b+S(a_0;1)+d+(H\cap\Gamma)\subset X$.
We have shown that
$$b+S(a_0,\a;1)+(H\cap\Gamma)=b+S(a_0;1)+S(\a;1)+(H\cap\Gamma)\subset X,$$
as desired.
\qed

\begin{corollary}
\label{congruence-free}
Suppose $X\subset G$ is a subvariety such that $X\cap\Gamma$ is an $F$-set.
Then $X\cap\Gamma$ can be expressed as an $F$-set where the basic groupless $F$-sets that appear do not involve congruence conditions, and the groups that appear are intersections of algebraic subgroups of $G$ over $\FF_q$ with $\Gamma$.
\end{corollary}

\proof
Suppose $\displaystyle X\cap\Gamma=\bigcup_{i=1}^l b_i+S_i(\a_i;\overline\delta_i)+H_{\Gamma}^i$.
For each $i$, let $H^i$ be the Zariski closure of $H_{\Gamma}^i$ in $G$.
Note that $H^i$ is over $\FF_q$ as it is over $\FF_q^{\alg}$ (being an algebraic subgroup of $G$) and over $K$ (having a Zariski dense set of points in $\Gamma\subset G(K)$).
By Lemma \ref{congruence-free-zc}, for each $i\leq l$, we have that, $b_i+S_i(\a_i;1)+(H^i\cap\Gamma)\subset X$.
On the other hand, it is clear that
$$b_i+S_i(\a_i;\overline\delta_i)+H_{\Gamma}^i\subset b_i+S_i(\a_i;1)+(H^i\cap\Gamma).$$
Hence, $\displaystyle X\cap\Gamma=\bigcup_{i=1}^l b_i+S_i(\a_i;1)+(H^i\cap\Gamma)$, as desired.
\qed

\smallskip

Next we deal with $F$-orbits of $F$-sets.
If $U$ is an $F$-set, then by the {\em $F$-orbit of $U$}, we mean the set $\displaystyle \bigcup_{n\geq 0}F^n U$.
Note that $F$-sets are not themselves closed under taking $F$-orbits.
Nevertheless, we have:

\begin{lemma}
\label{f-orbit-zc}
Suppose $U=b+S(\a;\overline\delta)+H_{\Gamma}$ is a translate of a basic $F$-set in $\Gamma$.
Let $H$ be the Zariski closure of $H_{\Gamma}$.
Then
$$S(b,\a;1)+(H\cap\Gamma)\subset\overline{\bigcup_{n\geq 0}F^n U}.$$
\end{lemma}

\proof
Note that for all $i\leq 0$
\begin{eqnarray*}
F^ib+S(\a;1)+(H\cap\Gamma)
& \subset & F^ib+\overline{S(\a;1)+(H\cap\Gamma)} \\
& = & F^ib+F^i[\overline{S(\a;1)+(H\cap\Gamma)}] \\
& = & F^i[b+\overline{S(\a;1)+(H\cap\Gamma)}] \\
& = & F^i(b+\overline{S(\a;\overline\delta)+H_{\Gamma}}) \\
& \subset & \overline{\bigcup_{n\geq 0}F^n U}
\end{eqnarray*}
where the first equality is because $\overline{S(\a;1)+H_{\Gamma}}$ is defined over $\FF_q$ by Lemma \ref{zc-basic-f-set}; and the third equality is by Lemma \ref{congruence-free-zc}.
Hence, $\displaystyle S(b,\a;1)+(H\cap\Gamma)\subset\overline{\bigcup_{n\geq 0}F^n U}$, as desired.
\qed

\smallskip

\begin{corollary}
\label{f-orbit-free}
Suppose $X\subset G$ is a subvariety such that $X\cap\Gamma$ is a finite union of $F$-sets and $F$-orbits of $F$-sets.
Then $X\cap\Gamma$ is in fact an $F$-set.
\end{corollary}

\proof
Suppose $\Sigma$ is the $F$-orbit of an $F$-set $U$, that appear in $X\cap\Gamma$.
We wish to replace $\Sigma$ by an $F$-set.
Taking finite unions, we may assume that $U$ is of the form $b+S(\a;\overline\delta)+H_{\Gamma}$.
By Lemma \ref{f-orbit-zc}, $S(b,\a;1)+(H\cap\Gamma)\subset\overline\Sigma\subset X$.
On the other hand, it is clear that $\Sigma\subset S(b,\a;1)+(H\cap\Gamma)$.
Hence we may replace $\Sigma$ with $S(b,\a;1)+(H\cap\Gamma)$ in the expression of $X\cap\Gamma$.
\qed

\medskip

Having understood some particular properties of $F$-sets in $\Gamma$, we now work toward showing that $X\cap\Gamma$ is indeed an $F$-set for every subvariety $X\subset G$.
First a preliminary lemma:

\begin{lemma}
\label{purity}
\begin{itemize}
\item[1.] The $F$-pure hull of $\Gamma$ in $G(K)$ -- i.e., the set of $g\in G(K)$ such that $F^ng\in\Gamma$ for some $n$ -- is a finitely generated group.
\item[2.] For all $n\geq 0$, $\Gamma/F^n\Gamma$ is finite.
\item[3.] For some $m>0$, $\Gamma\setminus F\Gamma\subset G(K)\setminus G(K^{q^m})$.
\end{itemize}
\end{lemma}

\proof
Part $1$ is clear in the case when $G$ is an abelian variety, since in that case $G(K)$ is itself a finitely generated group.
For the general case, we can choose $\R\subset K$ an integrally closed ring that is finitely genrated over $\FF_q$ (as a ring), and such that $\Gamma\leq G(\R)$.
Now $G(\R)$ is a finitely generated group, and since $K^q\cap\R=\R^q$, it is $F$-pure in $G(K)$.
Hence the $F$-pure hull of $\Gamma$ in $G(K)$, $\Gamma^\prime$, is a subgroup of $G(\R)$.
It follows that $\Gamma^\prime$ is itself a finitely generated group.

Recall that the multiplication by $q^n$ map on $G(K)$ is equal to $F^n\circ V^n$, where $V\colon G(K)\to G(K)$ is the Verschiebung map.
As $\Gamma^\prime$ is $F$-pure in $G(K)$, it follows that $V$ restricts to an endomorphism of $\Gamma^\prime$, and that $q^n\Gamma^\prime\subset F^n\Gamma^\prime$.
Now $\Gamma^\prime/q^n\Gamma^\prime$, being a finitely generated $\ZZ/q^n\ZZ$-module, is finite.
Hence $\Gamma^\prime/F^n\Gamma^\prime$ is finite for all $n$.
Now consider $\Gamma$ itself, and let $n$ be arbitrary.
As $\Gamma^\prime$ is a finitely generated group, for some $N>0$, $F^N\Gamma^\prime\subset\Gamma$.
So $\Gamma^\prime/F^n\Gamma$ is a quotient of the finite group $\Gamma^\prime/F^{n+N}\Gamma^\prime$.
Hence $\Gamma/F^n\Gamma$ is finite, establishing $2$.

Finally, for $3$, let $N>0$ again be such that $F^N\Gamma^\prime\subset\Gamma$, and let $m=N+1$.
Then $\Gamma\cap G(K^{q^m})\subset F^m\Gamma^\prime\subset F\Gamma$.
It follows that $\Gamma\setminus F\Gamma\subset G(K)\setminus G(K^{q^m})$.
\qed

\smallskip

\begin{remark}
\label{defzc}
In what follows we will implicitly use the following fact:
{\em 
Let $L$ be a field possible equipped with additional structure.
Suppose $Y$ is a variety defined over $L$, $\Upsilon \subseteq Y(L)$ is a set definable in L, and $\{X_b\}_{b \in B}$ is an algebraic family of subvarieties of $Y$ defined over $L$.
Then the condition that $X_b \cap \Upsilon$ is Zariski dense in $X_b$ is a type-definable 
condition on $b$.
}
Indeed, if $X_b \cap \Upsilon$ is {\em not} Zariski dense in $X_b$, then its Zariski closure $V\subsetneq X_b$, is over $L$, and $V(L)$ contains $X_b \cap \Upsilon$.
The existence of such a $V$ may be expressed as a countable disjuntion of formulas, hence the negation is type definable.
Indeed, that $V$ is contained in $X_b$ can be expressed by saying that a given set of polynomials defining $X_b$ is contained in the radical of the ideal generated by a given set of polynomials defining $V$.
That the containment is proper is equivalent to there being some natural number $d$, and an extension $L'/L$ of degree $d$, for which $V(L') \neq X_b(L')$.
As the extensions of degree $d$ may be uniformly encoded in $L$, this too can be expressed using countable many formulas.
\end{remark}

\smallskip

The following proposition was shown in \cite{abramovich-voloch} using a Hilbert scheme argument\footnote{This was their ``warm-up case''.}; and also follows from Hrushovski's proof of the function field Mordell-Lang conjecture in positive characteristic \cite{hrushovski}.
We present an elementary model-theoretic argument based on one that appeared in unpublished notes of the first author.

\begin{proposition}
\label{rel-ml}
Suppose $X\subset G$ is a subvariety such that $X\cap\Gamma$ is Zariski dense in $X$.
Then for some $\gamma\in G(K^{\alg})$, $\gamma+X$ is defined over $\FF_q^{\alg}$.
\end{proposition}

\proof
Note that $X$ is defined over $K$, since it has a dense intersection with $G(K)$.
Let $L$ denote the separable closure of $K$.
Since $K$ is a finitely generated extension of $\FF_q$, we have that $\displaystyle \bigcap_{n\geq 0}L^{q^n}=\FF_q^{\alg}$.

Let ${^*}L$ be an $\omega_1$-saturated elementary extension of $L$, and let
$\displaystyle k= \bigcap_{n \geq 0} ({^*}L)^{q^n}$.
Note that $L$ is linearly disjoint from $k$ over $\FF_q^{\alg}$.
Hence, $L^{\alg}=K^{\alg}$ is algebraically disjoint from $k$ over $\FF_q^{\alg}$.
It follows that
$(K^{\alg},\FF_q^{\alg})\preceq(({^*}L)^{\alg},k)$.

Fix a natural number $n$.
As $\Gamma$ is a finite union of cosets of $F^n\Gamma$, and $X\cap \Gamma$ is Zariski dense in $X$; for some $\gamma_n \in \Gamma$, $(\gamma_n+X)\cap F^n\Gamma$ is Zariski dense in $\gamma_n+X$.
Hence $(\gamma_n+X)\cap G(L^{q^n})$ is Zariski dense in $\gamma_n+X$.
As this is a type-definable property of $\gamma_n$, and by saturation of ${^*}L$, there is ${^*}\gamma\in G({^*}L)$ such that $({^*}\gamma+X)\cap G({^*}L^{q^n})$ is Zariski dense in ${^*}\gamma+X$, for all $n$.
By saturation again, $({^*}\gamma+X)\cap G(k)$ is Zariski dense in ${^*}\gamma+X$.
It follows that ${^*}\gamma+X$ is defined over $k$.

In particular, we have shown that some translate of $X$ by an element of $G(({^*}L)^{\alg})$ is defined over $k$.
This can be witnessed by a first order sentence with parameters from $L$, true of $(({^*}L)^{\alg},k)$.
By tranfer, it is true of $(K^{\alg},\FF_q^{\alg})$.
\qed

\medskip

We are now ready to prove the main result of this section: a version of the absolute Mordell-Lang conjecture for finitely generated $\ZZ[F]$-submodules of semiabelian varieties over finite fields.
This proof is based on an argument that was first presented by the second author in the e-print~\cite{Sca}.

\begin{theorem}
\label{ml}
If $X \subseteq G$ is any closed subvariety, then $X\cap\Gamma$ is an $F$-set.
The basic $F$-sets that appear have the form $S(a_1, \ldots, a_n;1) + (H\cap\Gamma)$ for some algebraic subgroup $H \leq G$ over $\FF_q$ and points $a_1, \ldots, a_n \in \Gamma$.
\end{theorem}

\begin{proof}
First of all, by Corollary \ref{congruence-free}, it is sufficient to show that $X\cap\Gamma$ is an $F$-set.
We proceed by induction on the dimension of $X$.
When $X$ is a finite set of points, the theorem is trivially true.

As the class of $F$-sets is closed under finite unions, we may assume that $X$ is irreducible.  
We may also assume that $X$ has a trivial stabilizer.
Let $H=\stab_G(X)$ be the stabilizer of $X$ in $G$ as an algebraic group.
Then $H$ is defined over $K$ and also ${\mathbb F}_q^{\alg}$ (as 
every algebraic subgroup of $G$ is defined over ${\mathbb F}_q^{\alg}$), and therefore over ${\mathbb F}_q$.
Let $\pi: G \to G/H =: \overline{G}$ be the quotient map.  Set $\overline{\Gamma} := \pi(\Gamma)$ and $\overline{X} := \pi(X)$.
Then, $\stab_{\overline{G}}(\overline{X})$ is trivial.
Assuming the result in the case of $\overline{X}$, we have that $\overline{X}\cap \overline{\Gamma}$ is an $F$-set.
Using the fact that the kernel of $\pi\upharpoonright_\Gamma$ stabilizes $X$, we have $X\cap \Gamma = \pi\upharpoonright_\Gamma^{-1}(\overline{X}\cap \overline{\Gamma})$.
So, $X\cap \Gamma$ is an $F$-set as long as the same is true of $\overline{X} \cap \overline{\Gamma}$.

Using Proposition \ref{rel-ml} we see that there is some $\gamma \in G(K^{\alg})$ such that $\gamma + X$ is defined over ${\mathbb F}_q^{\alg}$.
Let $K^\prime=K(\gamma)$.
Note that $\gamma + X$ is defined over $K^\prime\cap\FF_q^{\alg}=\FF_{q^r}$ for some $r>0$.
Let $\Gamma^\prime=\Gamma<\gamma>$, be the $\ZZ[F]$-submodule of $G(K^\prime)$ generated by $\Gamma$ and $\gamma$.
We view $\Gamma^\prime$ as an $F^r$-space.
Now $\Gamma$ is an $\ZZ[F^r]$-submodule of $\Gamma^\prime$; and hence an $F^r$-subspace of $\Gamma^\prime$ (using Proposition \ref{module=space}).
It follows that if $X\cap\Gamma^\prime$ is an $F^r$-set, then so is $X\cap\Gamma=(X\cap\Gamma^\prime)\cap\Gamma$.
Being contained in $\Gamma$, $X\cap\Gamma$ would then be an $F$-set in the sense of $\Gamma$ itself -- and we would be done.
That is, we may assume that $K^\prime=K$, $r=1$, and $\Gamma^\prime=\Gamma$.
In other words, we may assume that for some $\gamma\in\Gamma$, $\gamma+X$ is defined over $\FF_q$.
Translating by $-\gamma$, we may assume that $X$ itself is defined over ${\mathbb F}_q$.

We now claim that for some $N>0$, if $\xi\in\Gamma\setminus F\Gamma$, then $(\xi + X)\cap G(K^{q^N})$ is \emph{not} Zariski dense in $\xi +X$.
Suppose this were false, and let $(\xi_i)_{i\in\omega}$ be a sequence of points in $\Gamma\setminus F\Gamma$, and $(n_i)_{i\in\omega}$ a strictly increasing sequence of positive integers, such that $(\xi_i+X)\cap G(K^{q^{n_{i}}})$ is Zariski dense in $\xi_i+X$.
Using part $3$ of Lemma \ref{purity}, choose $m>0$ such that $\Gamma\setminus F\Gamma\subset G(K)\setminus G(K^{q^m})$.
Hence, each $\xi_i\in G(K)\setminus G(K^{q^m})$.
Let ${^*}K$ be an $\omega_1$-saturated elementary extension of $K$.
For each $i\in\omega$, we have that $\xi_i\in G({^*}K)\setminus G({^*}K^{q^m})$, and $(\xi_i+X)\cap G({^*}K^{q^{n_{i}}})$ is Zariski dense in $\xi_i+X$.
This being a type-definable property of $\xi_i$, and by saturation of ${^*}K$, we obtain ${^*}\xi\in G({^*}K)\setminus G({^*}K^{q^m})$ such that $({^*}\xi+X)\cap G({^*}K^{q^{n_{i}}})$ is Zariski dense in ${^*}\xi+X$, for all $i$.
By saturation again, $({^*}\xi+X)\cap G(k)$ is Zariski dense in ${^*}\xi+X$, where $\displaystyle k=\bigcap_{n\geq 0}{^*}K^{q^n}$.
It follows that ${^*}\xi+X$ is defined over $k$.
But as $X$ is defined over $\FF_q\subset k$, we obtain that for all $\sigma\in\aut(({^*}K^{\alg})/k)$, $\sigma({^*}\xi)-{^*}\xi$ stabilises $X$.
As $\stab X$ is trivial (even in $({^*}K)^{\alg}$), this means that $\sigma$ fixes ${^*}\xi$.
Hence, ${^*}\xi\in G(k)\subset G({^*}K^{q^m})$, which is a contradiction.

Now let $\xi\in\Gamma\setminus F\Gamma$, and choose coset representative $\eta_1,\dots \eta_l$ for $F^N\Gamma$ in $F\Gamma$.
Then each of $(\xi+\eta_i+X)\cap F^N\Gamma$ is not Zariski dense in $(\xi+\eta_i+X)$.
Let $Z_i$ be the Zariski closure of $(\xi+\eta_i+X)\cap F^N\Gamma$.
As $X$ is irreducible, this means that each $Z_i$ has dimension strictly less than that of $X$.
But
$$\dim\overline{(\xi+X)\cap F\Gamma}=\dim(\bigcup_{i=1}^l\overline{(\xi+\eta_i+X)\cap F^N\Gamma})=\dim(\bigcup_{i=1}^lZ_i)<\dim X.$$
We have shown that for all $\xi\in\Gamma\setminus F\Gamma$, the Zariski closure of $(\xi+X)\cap F\Gamma$ is of strictly smaller dimension than that of $X$.

We are now ready to finish the proof of the theorem.
Note that $\Gamma$ can be written as $\displaystyle (G({\mathbb F}_q) \cap \Gamma) \cup \bigcup_{n \geq 0} F^n(\Gamma \setminus F\Gamma)$.
Fix coset representatives $\gamma_1, \ldots, \gamma_\ell$ for the \emph{nonzero} cosets of $F \Gamma$ in $\Gamma$.
We have,
\begin{eqnarray*}
X\cap \Gamma
& = & X\cap[(G({\mathbb F}_q) \cap \Gamma) \cup \bigcup_{n \geq 0} F^n(\Gamma \setminus F\Gamma)] \\
& = & (X(\FF_q)\cap\Gamma) \cup \bigcup_{n \geq 0} F^n(X\cap (\Gamma \setminus F\Gamma)) \\
& = & (X(\FF_q)\cap\Gamma) \cup \bigcup_{i=1}^\ell \bigcup_{n \geq 0}F^n(\gamma_i + [(X - \gamma_i)\cap F\Gamma])
\end{eqnarray*}
where the second equality uses the fact that $X$ is defined over $\FF_q$, and hence is fixed by powers of $F$.
For each $i\leq\ell$, $\gamma_i\in\Gamma\setminus F\Gamma$.
We have already pointed out that in this case the dimension of the Zariski closure of $(X - \gamma_i)\cap F\Gamma$ is strictly less than that of $X$.
Hence by induction, $\gamma_i + [(X - \gamma_i)\cap F\Gamma]$ is an $F$-set for each $i\leq\ell$.
So the above equation displays $X\cap\Gamma$ as the union of a finite set with a finite union of $F$-orbits of $F$-sets.
By Corollary \ref{f-orbit-free}, we obtain that $X\cap\Gamma$ is an $F$-set.
\end{proof}

\smallskip

\begin{corollary}
\label{stable-predicate}
Suppose  $\U$ is an algebraically closed field extending $K$.  Then $\th(\U,+,\times,\Gamma)$ is stable.
\end{corollary}

\proof
Theorem \ref{ml} implies that the structure induced on $\Gamma$ by $(\U,+,\times,\Gamma)$ is a fragment of $(\Gamma,\S)$.
More precisely, denote by $\Gamma_{\ind}$ the structure whose universe is $\Gamma$ and whose relations are sets of the form $D\cap\Gamma^n$, where $D\subset G(\U)^n$ is definable with parameters in $(\U,+,\times)$.
Then every definable set in $\Gamma_{\ind}$ is definable in $(\Gamma,\S)$.
Indeed, by Theorem \ref{ml}, $D\cap\Gamma^n$ is an $F$-set when $D\subset G(\U)^n$ is a closed subvariety; and as every set $D\subset G(\U)^n$ definable in $(\U,+,\times)$ is a boolean combination of closed sets, $D\cap\Gamma^n$ is a boolean combination of $F$-sets.
By Theorem \ref{stable}, $\Gamma_{\ind}$ is stable.
By a result of Pillay's from \cite{pillay}, $(\U,+,\times,\Gamma)$ is stable.
\qed

\medskip

We also obtain a uniform version of Theorem \ref{ml}.

\begin{corollary}
\label{uniform-ml}
Suppose $\{ X_b \}_{b \in B}$ is an algebraic family of closed subvarieties of $G$.
Then there are basic groupless $F$-sets $T_1, \ldots, T_m \subseteq \Gamma$ and algebraic subgroups $H_1, \ldots, H_m \leq G$ over $\FF_q$, such that for any $b \in B$ there is an $I \subseteq \{1, \ldots, m \}$ and points $(\gamma_i)_{i\in I}$ from $\Gamma$, such that
$$X_b\cap \Gamma = \bigcup_{i \in I} \gamma_i + T_i + (H_i \cap \Gamma)$$.
\end{corollary}

\begin{proof}
Let $\U$ be an algebraically closed field extending $K$, over which the family $\{ X_b \}_{b \in B}$ is defined.
Let $({^*}\U,+,\times,{^*}\Gamma,{^*}\S)$ be a sufficiently saturated elementary extension of $(\U,+,\times,\Gamma,\S)$.
By saturation it suffices to prove that for every closed subvariety $X\subset G$ (over ${^*}\U$), $X\cap{^*}\Gamma$ is a finite union of translates of sets of the form $T+(H\cap{^*}\Gamma)$, where $T\in{^*}\S$ is the interpretation of a groupless $F$-set, and $H$ is an algebraic subgroup of $G$ over $\FF_q$.
As in Theorem \ref{ml}, we may assume that $X$ is irreducible, that $X\cap{^*}\Gamma$ is Zariski dense in $X$, and that the stabiliser of $X$ is trivial.
We will actually show that $X$ is a translate (by an element of ${^*}\Gamma$) of a subvariety defined over $\FF_q^{\alg}\subset\U$.
This is sufficient, as by Theorem \ref{ml} and transfer, the intersection of a subvariety over $\U$ with ${^*}\Gamma$ is of the required form.

We obtain an induced structure ${^*}\Gamma_{\ind}$ on ${^*}\Gamma$, which is an elementary extension $\Gamma_{\ind}$.
Namely, the relations are sets of the form $D\cap{^*}\Gamma^n$ where $D\subset G({^*}\U)^n$ is definable in $({^*}\U,+,\times)$ with paramters from $\U$.
Note that ${^*}\Gamma_{\ind}$ is a fragment of $({^*}\Gamma,{^*}\S)$.
Moreover, the stability of $\Gamma_{\ind}$ implies that every set in ${^*}\Gamma$ definable in $({^*}\U,+,\times,{^*}\Gamma)$ is already definable in ${^*}\Gamma_{\ind}$.\footnote{This follows from Pillay's arguments in \cite{pillay}; specifically from Lemma $2.9$ and the discussion surrounding it.}
Hence $X\cap{^*}\Gamma$ is a finite boolean combination of translates of sets from ${^*}\S$ by elements of ${^*}\Gamma$.
Possibly translating $X$ by an element of ${^*}\Gamma$, and using the irreducibility of $X$, we have that there are basic groupless $S,T_1,\dots T_n\in{^*}\S$, $F$-subspaces $H,L_1,\dots L_n\in{^*}\S$, and points $a_1,\dots a_n\in{^*}\Gamma$, such that $\displaystyle (S+H)\setminus(\bigcup_{i=1}^n(a_i+T_i+L_i))$ is a Zariski dense subset of $X$.
Let $U=S+H$ and $\displaystyle V=\bigcup_{i=1}^n(a_i+T_i+L_i)$.
Intersecting, we may assume that $V\subset U$ and that each $L_i$ is a subgroup of $H$.
We begin with some reductions.

{\em We may assume that each $L_i$ has finite index in $H$.}
Let $V'\subset V$ be the sub-union involoving only those $L_i$ which are of infinite index in $H$.
By Lemma \ref{finind}, $H$ is not contained in $V'-V'$.
Hence, for some $h\in H$, $V'\cap(h+V')$ is empty.
It follows that
$$(U\setminus V)\cup(h+[U\setminus V])=U\setminus W,$$
for some $W$ a finite union of translates of $F$-sets,  where the subgroups involved are of finite index in $H$.
Suppose we knew the result for $U\setminus W$; that is, its Zariski closure is a finite union of translates of subvarieties over $\FF_q^{\alg}$.
Taking Zariski closures of both sides, we would have that this is true of $X\cup(h+X)$.
As $X$ is irreducible, it would itself be a translate of a subvariety over $\FF_q^{\alg}$, as desired.

{\em The groups $H$ and $L_i$ are finite.}
Let $L$ be the intersection of the $L_i$'s.
Then $L$ is of finite index in $H$ and in each $L_i$.
Moreover, $L$ stabilises $U$ and $V$, and hence $X$.
Since $X$ has trivial stabiliser, $L=0$, and $H$ as well as the $L_i$'s must be finite.

It follows that $U$ and $V$ are finite unions of translates of basic groupless $F$-sets.
We have thus reduced the corollary to the following lemma (where we are renaming the symbols $S$ and $T$):

\begin{lemma}
If $T\subset S$ are finite unions of translates of basic groupless $F$-sets in ${^*}\Gamma$, then the Zariski closure of $S\setminus T$ is a finite union of translates of subvarieties defined over $\FF_q^{\alg}$.
\end{lemma}

\proof
We fix $\delta$ such that the basic groupless $F$-sets involved in $S$ are from $O(\delta)$, and we proceed by induction on the maximal $\delta$-dimension of the basic groupless $F$-sets involved in $S$.
The case when this dimension is $0$ is trivial.

We begin with a few reductions.
Suppose $Z$ is an irreducible component of the Zariski closure of $S\setminus T$.
Then $(S\setminus T)\cap Z$ is Zariski dense in $Z$.
As $Z\cap{^*}\Gamma$ is a boolean combination of translates of basic $F$-sets, $(S\setminus T)\cap Z=(S\setminus T)\cap Z\cap{^*}\Gamma$ is again of the form $S'\setminus T'$, where $S'$ and $T'$ are still unions of translates of groupless $F$-sets, and $S'\subset S$.
Hence we may assume that the Zariski closure of $S\setminus T$, $Z$, is irreducible.
Moreover, taking finite unions and translating, we may assume that $S$ is a basic groupless $F$-set of $\delta$-degree $1$.
Finally, by induction, and Proposition~\ref{closure-congruence-free}, we may assume that the basic groupless $F$-sets involved in $T$ have $\delta$-dimension strictly less than that of $S$.

Now, $\displaystyle F(S\setminus T)=FS\setminus FT=[FS\setminus(T\cup FT)]\cup[T\setminus FT]$.
Taking Zariski closures, and using the fact that $Z$ is irreducible, we have that $FZ$ is either the Zariski closure of $T\setminus FT$, or of $FS\setminus(T\cap FT)$.
In the former case, by induction, $FZ$, and hence $Z$ is a translate of a subvariety defined over $\FF_q^{\alg}$.
In the latter case, we have that
$$FZ=\overline{FS\setminus(T\cup FT)}\subset\overline{S\setminus T}=Z,$$
where we are using the fact that $S$ is a basic groupless $F$-set, and hence $FS\subset S$.
It follows that $Z$ is fixed by $F$, and hence is over $\FF_q$.

This completes the proof of the lemma and Corollary \ref{uniform-ml}.
\qed

\end{proof}

\end{document}